\magnification=\magstep1
\input amstex
\def\PaperSize{letter}          
 

\documentstyle{amsppt}


\brokenpenalty=10000
\clubpenalty=10000
\widowpenalty=10000

\catcode`\@=11

\def\keyboarder#1{}%

\def\pagewidth#1{\hsize#1
  \captionwidth@24pc}

\pagewidth{30pc}

\parindent=18\p@
\normalparindent\parindent

\parskip=\z@

\def\foliofont@{\sevenrm}
\def\headlinefont@{\sevenpoint}


\def\specialheadfont@{\elevenpoint\smc}
\def\headfont@{\bf}
\def\subheadfont@{\bf}
\def\refsheadfont@{\bf}
\def\abstractfont@{\smc}
\def\proclaimheadfont@{\smc}
\def\xcaheadfont@{\smc}
\def\captionfont@{\smc}
\def\citefont@{\bf}
\def\refsfont@{\eightpoint}

\font\sixsy=cmsy6


\font@\twelverm=cmr10 scaled \magstep1
\font@\twelvebf=cmbx10 scaled \magstep1
\font@\twelveit=cmti10 scaled \magstep1
\font@\twelvesl=cmsl10 scaled \magstep1
\font@\twelvesmc=cmcsc10 scaled \magstep1
\font@\twelvett=cmtt10 scaled \magstep1
\font@\twelvei=cmmi10 scaled \magstep1
\font@\twelvesy=cmsy10 scaled \magstep1
\font@\twelveex=cmex10 scaled \magstep1
\font@\twelvemsb=msbm10 scaled \magstep1
\font@\twelveeufm=eufm10 scaled \magstep1

\newtoks\twelvepoint@
\def\twelvepoint{\normalbaselineskip14\p@
  \abovedisplayskip12\p@ plus3\p@ minus9\p@
  \belowdisplayskip\abovedisplayskip
  \abovedisplayshortskip\z@ plus3\p@
  \belowdisplayshortskip7\p@ plus3\p@ minus4\p@
  \textonlyfont@\rm\twelverm \textonlyfont@\it\twelveit
  \textonlyfont@\sl\twelvesl \textonlyfont@\bf\twelvebf
  \textonlyfont@\smc\twelvesmc \textonlyfont@\tt\twelvett
  \ifsyntax@ \def\big##1{{\hbox{$\left##1\right.$}}}%
    \let\Big\big \let\bigg\big \let\Bigg\big
  \else
    \textfont\z@\twelverm  \scriptfont\z@\eightrm
       \scriptscriptfont\z@\sixrm
    \textfont\@ne\twelvei  \scriptfont\@ne\eighti
       \scriptscriptfont\@ne\sixi
    \textfont\tw@\twelvesy \scriptfont\tw@\eightsy
       \scriptscriptfont\tw@\sixsy
    \textfont\thr@@\twelveex \scriptfont\thr@@\eightex
        \scriptscriptfont\thr@@\eightex
    \textfont\itfam\twelveit \scriptfont\itfam\eightit
        \scriptscriptfont\itfam\eightit
    \textfont\bffam\twelvebf \scriptfont\bffam\eightbf
        \scriptscriptfont\bffam\sixbf
    \textfont\msbfam\twelvemsb \scriptfont\msbfam\eightmsb
        \scriptscriptfont\msbfam\sixmsb
    \textfont\eufmfam\twelveeufm \scriptfont\eufmfam\eighteufm
        \scriptscriptfont\eufmfam\sixeufm
    \setbox\strutbox\hbox{\vrule height8.5\p@ depth3.5\p@ width\z@}%
    \setbox\strutbox@\hbox{\lower.5\normallineskiplimit\vbox{%
        \kern-\normallineskiplimit\copy\strutbox}}%
    \setbox\z@\vbox{\hbox{$($}\kern\z@}\bigsize@1.2\ht\z@
  \fi
  \normalbaselines\rm\dotsspace@1.5mu\ex@.2326ex\jot3\ex@
  \the\twelvepoint@}

\font@\elevenrm=cmr10 scaled \magstephalf
\font@\elevenbf=cmbx10 scaled \magstephalf
\font@\elevenit=cmti10 scaled \magstephalf
\font@\elevensl=cmsl10 scaled \magstephalf
\font@\elevensmc=cmcsc10 scaled \magstephalf
\font@\eleventt=cmtt10 scaled \magstephalf
\font@\eleveni=cmmi10 scaled \magstephalf
\font@\elevensy=cmsy10 scaled \magstephalf
\font@\elevenex=cmex10 scaled \magstephalf
\font@\elevenmsb=msbm10 scaled \magstephalf
\font@\eleveneufm=eufm10 scaled \magstephalf

\newtoks\elevenpoint@
\def\elevenpoint{\normalbaselineskip13\p@
  \abovedisplayskip12\p@ plus3\p@ minus9\p@
  \belowdisplayskip\abovedisplayskip
  \abovedisplayshortskip\z@ plus3\p@
  \belowdisplayshortskip7\p@ plus3\p@ minus4\p@
  \textonlyfont@\rm\elevenrm \textonlyfont@\it\elevenit
  \textonlyfont@\sl\elevensl \textonlyfont@\bf\elevenbf
  \textonlyfont@\smc\elevensmc \textonlyfont@\tt\eleventt
  \ifsyntax@ \def\big##1{{\hbox{$\left##1\right.$}}}%
    \let\Big\big \let\bigg\big \let\Bigg\big
  \else
    \textfont\z@\elevenrm  \scriptfont\z@\eightrm
       \scriptscriptfont\z@\sixrm
    \textfont\@ne\eleveni  \scriptfont\@ne\eighti
       \scriptscriptfont\@ne\sixi
    \textfont\tw@\elevensy \scriptfont\tw@\eightsy
       \scriptscriptfont\tw@\sixsy
    \textfont\thr@@\elevenex \scriptfont\thr@@\eightex
        \scriptscriptfont\thr@@\eightex
    \textfont\itfam\elevenit \scriptfont\itfam\eightit
        \scriptscriptfont\itfam\eightit
    \textfont\bffam\elevenbf \scriptfont\bffam\eightbf
        \scriptscriptfont\bffam\sixbf
    \textfont\msbfam\elevenmsb \scriptfont\msbfam\eightmsb
        \scriptscriptfont\msbfam\sixmsb
    \textfont\eufmfam\eleveneufm \scriptfont\eufmfam\eighteufm
        \scriptscriptfont\eufmfam\sixeufm
    \setbox\strutbox\hbox{\vrule height8.5\p@ depth3.5\p@ width\z@}%
    \setbox\strutbox@\hbox{\lower.5\normallineskiplimit\vbox{%
        \kern-\normallineskiplimit\copy\strutbox}}%
    \setbox\z@\vbox{\hbox{$($}\kern\z@}\bigsize@1.2\ht\z@
  \fi
  \normalbaselines\rm\dotsspace@1.5mu\ex@.2326ex\jot3\ex@
  \the\elevenpoint@}

\addto\tenpoint{\normalbaselineskip12\p@
 \abovedisplayskip6\p@ plus6\p@ minus0\p@
 \belowdisplayskip6\p@ plus6\p@ minus0\p@
 \abovedisplayshortskip0\p@ plus3\p@ minus0\p@
 \belowdisplayshortskip2\p@ plus3\p@ minus0\p@
 \ifsyntax@
 \else
  \setbox\strutbox\hbox{\vrule height9\p@ depth4\p@ width\z@}%
  \setbox\strutbox@\hbox{\vrule height8\p@ depth3\p@ width\z@}%
 \fi
 \normalbaselines\rm}

\newtoks\sevenpoint@
\def\sevenpoint{\normalbaselineskip9\p@
 \textonlyfont@\rm\sevenrm \textonlyfont@\it\sevenit
 \textonlyfont@\sl\sevensl \textonlyfont@\bf\sevenbf
 \textonlyfont@\smc\sevensmc \textonlyfont@\tt\seventt
  \textfont\z@\sevenrm \scriptfont\z@\sixrm
       \scriptscriptfont\z@\fiverm
  \textfont\@ne\seveni \scriptfont\@ne\sixi
       \scriptscriptfont\@ne\fivei
  \textfont\tw@\sevensy \scriptfont\tw@\sixsy
       \scriptscriptfont\tw@\fivesy
  \textfont\thr@@\sevenex \scriptfont\thr@@\sevenex
   \scriptscriptfont\thr@@\sevenex
  \textfont\itfam\sevenit \scriptfont\itfam\sevenit
   \scriptscriptfont\itfam\sevenit
  \textfont\bffam\sevenbf \scriptfont\bffam\sixbf
   \scriptscriptfont\bffam\fivebf
  \textfont\msbfam\sevenmsb \scriptfont\msbfam\sixmsb
   \scriptscriptfont\msbfam\fivemsb
  \textfont\eufmfam\seveneufm \scriptfont\eufmfam\sixeufm
   \scriptscriptfont\eufmfam\fiveeufm
 \setbox\strutbox\hbox{\vrule height7\p@ depth3\p@ width\z@}%
 \setbox\strutbox@\hbox{\raise.5\normallineskiplimit\vbox{%
   \kern-\normallineskiplimit\copy\strutbox}}%
 \setbox\z@\vbox{\hbox{$($}\kern\z@}\bigsize@1.2\ht\z@
 \normalbaselines\sevenrm\dotsspace@1.5mu\ex@.2326ex\jot3\ex@
 \the\sevenpoint@}

\newskip\abovespecheadskip   \abovespecheadskip20\p@ plus8\p@ minus2\p@
\newdimen\belowspecheadskip  \belowspecheadskip6\p@
\outer\def\specialhead{%
  \add@missing\endroster \add@missing\enddefinition
  \add@missing\enddemo \add@missing\endexample
  \add@missing\endproclaim
  \penaltyandskip@{-200}\abovespecheadskip
  \begingroup\interlinepenalty\@M\rightskip\z@ plus\hsize
  \let\\\linebreak
  \specialheadfont@\raggedcenter@\noindent}

\let\varindent@\indent

\let\subsubhead\relax
\outer\def\subsubhead{%
  \add@missing\endroster \add@missing\enddefinition
  \add@missing\enddemo
  \add@missing\endexample \add@missing\endproclaim
  \let\savedef@\subsubhead \let\subsubhead\relax
  \def\subsubhead##1\endsubsubhead{\restoredef@\subsubhead
      {\def\usualspace{\/{\subsubheadfont@\enspace}}%
    \subsubheadfont@##1\unskip\frills@{.\enspace}}\ignorespaces}%
  \nofrillscheck\subsubhead}


\newskip\abstractindent 	\abstractindent=3pc
\long\def\block #1\endblock{\vskip 6pt
	{\leftskip=\abstractindent \rightskip=\abstractindent
	\noindent #1\endgraf}\vskip 6pt}

\long\def\ext #1\endext{\removelastskip\block #1\endblock}

\outer\def\xca{\let\savedef@\xca \let\xca\relax
  \add@missing\endproclaim \add@missing\endroster
  \add@missing\endxca \envir@stack\endxca
  \def\xca##1{\restoredef@\xca
    \penaltyandskip@{-100}\medskipamount
    \bgroup{\def\usualspace{{\xcaheadfont@\enspace}}%
      \varindent@\xcaheadfont@\ignorespaces##1\unskip
      \frills@{.\xcaheadfont@\enspace}}%
      \ignorespaces}%
  \nofrillscheck\xca}
\def\endxca{\egroup\revert@envir\endxca
  \par\medskip}

\def\remarkheadfont@{\smc}
\def\remark{\let\savedef@\remark \let\remark\relax
  \add@missing\endroster \add@missing\endproclaim
  \envir@stack\endremark
  \def\remark##1{\restoredef@\remark
    \penaltyandskip@{-100}\medskipamount
    {\def\usualspace{{\remarkheadfont@\enspace}}%
     \varindent@\remarkheadfont@\ignorespaces##1\unskip
     \frills@{.\enspace}}\rm
    \ignorespaces}\nofrillscheck\remark}
\def\endremark{\par\revert@envir\endremark\medskip}

\def\qed{\ifhmode\unskip\nobreak\fi\hfill
  \ifmmode\square\else$\m@th\square$\fi}


\newdimen\rosteritemsep
\rosteritemsep=.5pc

\newdimen\rosteritemitemwd
\newdimen\rosteritemitemitemwd

\newbox\setwdbox
\setbox\setwdbox\hbox{0.}\rosteritemwd=\wd\setwdbox
\setbox\setwdbox\hbox{0.\hskip.5pc(c)}\rosteritemitemwd=\wd\setwdbox
\setbox\setwdbox\hbox{0.\hskip.5pc(c)\hskip.5pc(iii)}%
  \rosteritemitemitemwd=\wd\setwdbox

\def\roster{%
  \envir@stack\endroster
  \edef\leftskip@{\leftskip\the\leftskip}%
  \relaxnext@
  \rostercount@\z@
  \def\item{\FN@\rosteritem@}%
  \def\itemitem{\FN@\rosteritemitem@}%
  \def\itemitemitem{\FN@\rosteritemitemitem@}%
  \DN@{\ifx\next\runinitem\let\next@\nextii@
    \else\let\next@\nextiii@
    \fi\next@}%
  \DNii@\runinitem
    {\unskip
     \DN@{\ifx\next[\let\next@\nextii@
       \else\ifx\next"\let\next@\nextiii@\else\let\next@\nextiv@\fi
       \fi\next@}%
     \DNii@[####1]{\rostercount@####1\relax
       \therosteritem{\number\rostercount@}~\ignorespaces}%
     \def\nextiii@"####1"{{\rm####1}~\ignorespaces}%
     \def\nextiv@{\therosteritem1\rostercount@\@ne~}%
     \par@\firstitem@false
     \FN@\next@}
  \def\nextiii@{\par\par@
    \penalty\@m\vskip-\parskip
    \firstitem@true}%
  \FN@\next@}

\def\rosteritem@{\iffirstitem@\firstitem@false
  \else\par\vskip-\parskip
  \fi
  \leftskip\rosteritemwd \advance\leftskip\normalparindent
  \advance\leftskip.5pc \noindent
  \DNii@[##1]{\rostercount@##1\relax\itembox@}%
  \def\nextiii@"##1"{\def\therosteritem@{\rm##1}\itembox@}%
  \def\nextiv@{\advance\rostercount@\@ne\itembox@}%
  \def\therosteritem@{\therosteritem{\number\rostercount@}}%
  \ifx\next[\let\next@\nextii@
  \else\ifx\next"\let\next@\nextiii@\else\let\next@\nextiv@\fi
  \fi\next@}

\def\itembox@{\llap{\hbox to\rosteritemwd{\hss
  \kern\z@ 
  \therosteritem@}\hskip.5pc}\ignorespaces}

\def\therosteritem#1{\rom{\ignorespaces#1.\unskip}}

\def\rosteritemitem@{\iffirstitem@\firstitem@false
  \else\par\vskip-\parskip
  \fi
  \leftskip\rosteritemitemwd \advance\leftskip\normalparindent
  \advance\leftskip.5pc \noindent
  \DNii@[##1]{\rostercount@##1\relax\itemitembox@}%
  \def\nextiii@"##1"{\def\therosteritemitem@{\rm##1}\itemitembox@}%
  \def\nextiv@{\advance\rostercount@\@ne\itemitembox@}%
  \def\therosteritemitem@{\therosteritemitem{\number\rostercount@}}%
  \ifx\next[\let\next@\nextii@
  \else\ifx\next"\let\next@\nextiii@\else\let\next@\nextiv@\fi
  \fi\next@}

\def\itemitembox@{\llap{\hbox to\rosteritemitemwd{\hss
  \kern\z@ 
  \therosteritemitem@}\hskip.5pc}\ignorespaces}

\def\therosteritemitem#1{\rom{(\ignorespaces#1\unskip)}}

\def\rosteritemitemitem@{\iffirstitem@\firstitem@false
  \else\par\vskip-\parskip
  \fi
  \leftskip\rosteritemitemitemwd \advance\leftskip\normalparindent
  \advance\leftskip.5pc \noindent
  \DNii@[##1]{\rostercount@##1\relax\itemitemitembox@}%
  \def\nextiii@"##1"{\def\therosteritemitemitem@{\rm##1}\itemitemitembox@}%
  \def\nextiv@{\advance\rostercount@\@ne\itemitemitembox@}%
  \def\therosteritemitemitem@{\therosteritemitemitem{\number\rostercount@}}%
  \ifx\next[\let\next@\nextii@
  \else\ifx\next"\let\next@\nextiii@\else\let\next@\nextiv@\fi
  \fi\next@}

\def\itemitemitembox@{\llap{\hbox to\rosteritemitemitemwd{\hss
  \kern\z@ 
  \therosteritemitemitem@}\hskip.5pc}\ignorespaces}

\def\therosteritemitemitem#1{\rom{(\ignorespaces#1\unskip)}}

\def\endroster{\relaxnext@
  \revert@envir\endroster 
  \par\leftskip@
  \penalty-50 
  \DN@{\ifx\next\Runinitem\let\next@\relax
    \else\nextRunin@false\let\item\plainitem@
      \ifx\next\par
        \DN@\par{\everypar\expandafter{\the\everypartoks@}}%
      \else
        \DN@{\noindent\everypar\expandafter{\the\everypartoks@}}%
      \fi
    \fi\next@}%
  \FN@\next@}


\def\address#1\endaddress{\global\advance\addresscount@\@ne
  \expandafter\gdef\csname address\number\addresscount@\endcsname
  {\vskip12\p@ minus6\p@\indent\addressfont@\smc\ignorespaces#1\par}}


\def\curraddr{\let\savedef@\curraddr
  \def\curraddr##1\endcurraddr{\let\curraddr\savedef@
  \toks@\expandafter\expandafter\expandafter{%
       \csname address\number\addresscount@\endcsname}%
  \toks@@{##1}%
  \expandafter\xdef\csname address\number\addresscount@\endcsname
  {\the\toks@\endgraf\noexpand\nobreak
    \indent\noexpand\addressfont@{\noexpand\rm
    \frills@{{\noexpand\it Current address\noexpand\/}:\space}%
    \def\noexpand\usualspace{\space}\the\toks@@\unskip}}}%
  \nofrillscheck\curraddr}

\def\email{\let\savedef@\email
  \def\email##1\endemail{\let\email\savedef@
  \toks@{\def\usualspace{{\it\enspace}}\endgraf\indent\addressfont@}%
  \toks@@{{\tt ##1}\par}%
  \expandafter\xdef\csname email\number\addresscount@\endcsname
  {\the\toks@\frills@{{\noexpand\it E-mail address\noexpand\/}:%
     \noexpand\enspace}\the\toks@@}}%
  \nofrillscheck\email}

\def\rom#1{{\rm #1}}

\def\bysame{\by\hbox to2pc{\hrulefill}\thinspace\kern\z@}

\def\refstyle#1{\uppercase{%
  \gdef\refstyle@{#1}%
  \if#1A\relax \def\keyformat##1{[##1]\enspace\hfil}%
  \else\if#1B\relax
    \refindentwd2pc
    \def\keyformat##1{\aftergroup\kern
              \aftergroup-\aftergroup\refindentwd}%
  \else\if#1C\relax
    \def\keyformat##1{\hfil##1.\enspace}%
  \fi\fi\fi}
}

\refstyle{A}

 \relax

\catcode`\@=11

\def\pretitle{\null\vskip74pt}

\def\addressfont@{\eightpoint}


%
%

\define\issueinfo#1#2#3#4{%
  \def\issuevol@{#1}\def\issueno@{#2}%
  \def\issuemonth@{#3}\def\issueyear@{#4}}

\define\originfo#1#2#3#4{\def\origvol@{#1}\def\origno@{#2}%
  \def\origmonth@{#3}\def\origyear@{#4}}

\define\copyrightinfo#1#2{\def\copyrightyear{#1}\def\crholder@{#2}}

\define\pagespan#1#2{\pageno=#1\def\start@page{#1}\def\end@page{#2}}

\issueinfo{00}{0}{}{1997}
\originfo{00}{0}{}{1997}
\copyrightinfo{\issueyear@}{American Mathematical Society}
\pagespan{000}{000}
\pageno=209 


\def\nojourlogo{\let\jourlogo\empty@}

\def\journame{AMS Proceedings Style}
\def\volinfo{Volume {\sixbf\issuevol@}, \issueyear@}
\let\jourlogoextra@\empty@
\let\jourlogoright@\empty@

\def\jourlogofont@{\sixrm\baselineskip7\p@\relax}
\def\jourlogo{%
  \vbox to\headlineheight{%
    \parshape\z@ \leftskip\z@ \rightskip\z@
    \parfillskip\z@ plus1fil\relax
    \jourlogofont@ \frenchspacing
    \line{\vtop{\parindent\z@ \hsize=.5\hsize
      \journame\newline\volinfo\jourlogoextra@
	\newline Revised November 2001\endgraf\vss}%
      \hfil
      \jourlogoright@
    }%
    \vss}%
}

\def\issn#1{\gdef\theissn{#1}}
\issn{0000-0000}


\def\copyrightline@{%
  \rightline{\sixrm \textfont2=\sixsy \copyright\copyrightyear\ \crholder@}}

\def\logo@{\copyrightline@}

\def\titlefont{%
 \ifsyntax@\else \twelvepoint\bf \fi }

\def\authorfont{%
  \ifsyntax@
  \else \elevenpoint
  \fi}

\def\title{\let\savedef@\title
  \def\title##1\endtitle{\let\title\savedef@
    \global\setbox\titlebox@\vtop{\titlefont\bf
      \raggedcenter@\frills@{##1}\endgraf}%
    \ifmonograph@
      \edef\next{\the\leftheadtoks}\ifx\next\empty \leftheadtext{##1}\fi
    \fi
    \edef\next{\the\rightheadtoks}\ifx\next\empty \rightheadtext{##1}\fi
  }%
  \nofrillscheck\title}

\def\author#1\endauthor{\global\setbox\authorbox@
  \vbox{\authorfont\raggedcenter@
    {\ignorespaces#1\endgraf}}\relaxnext@
  \edef\next{\the\leftheadtoks}%
  \ifx\next\empty\expandafter\uppercase{\leftheadtext{#1}}\fi}

\def\abstract{\let\savedef@\abstract
  \def\abstract{\let\abstract\savedef@
    \setbox\abstractbox@\vbox\bgroup\indenti=3pc\noindent$$\vbox\bgroup
      \def\envir@end{\endabstract}\advance\hsize-2\indenti
      \def\usualspace{\enspace}\eightpoint \noindent
      \frills@{{\abstractfont@ Abstract.\enspace}}}%
  \nofrillscheck\abstract}

\def\dedicatory #1\enddedicatory{\def\preabstract{{\vskip 20\p@
  \eightpoint\it \raggedcenter@#1\endgraf}}}

\outer\def\endtopmatter{\add@missing\endabstract
  \edef\next{\the\leftheadtoks}%
  \ifx\next\empty@
    \expandafter\leftheadtext\expandafter{\the\rightheadtoks}%
  \fi
  \ifx\thesubjclass@\empty@\else \makefootnote@{}{\thesubjclass@}\fi
  \ifx\thekeywords@\empty@\else \makefootnote@{}{\thekeywords@}\fi
  \ifx\thethanks@\empty@\else \makefootnote@{}{\thethanks@}\fi
  \inslogo@
  \pretitle
  \box\titlebox@
  \topskip10pt
  \preauthor
  \ifvoid\authorbox@\else \vskip16\p@ plus6\p@ minus0\p@\unvbox\authorbox@\fi
  \predate
  \ifx\thedate@\empty\else \vskip6\p@ plus2\p@ minus0\p@
    \line{\hfil\thedate@\hfil}\fi
  \setabstract@
  \nobreak
  \ifvoid\tocbox@\else\vskip1.5pc plus.5pc \unvbox\tocbox@\fi
  \prepaper
  \vskip36\p@\tenpoint
}

\def\setabstract@{%
  \preabstract
  \ifvoid\abstractbox@\else \vskip20\p@ \unvbox\abstractbox@ \fi
}



\begingroup
\let\specialhead\relax
\let\head\relax
\let\subhead\relax
\let\subsubhead\relax
\let\title\relax
\let\chapter\relax

\gdef\setwidest@#1#2{%
   \ifx#1\head\setbox\tocheadbox@\hbox{#2.\enspace}%
   \else\ifx#1\subhead\setbox\tocsubheadbox@\hbox{#2.\enspace}%
   \else\ifx#1\subsubhead\setbox\tocsubheadbox@\hbox{#2.\enspace}%
   \else\ifx#1\key
       \if C\refstyle@ \else\refstyle A\fi
       \setboxz@h{\refsfont@\keyformat{#2}}%
       \refindentwd\wd\z@
   \else\ifx#1\no\refstyle C%
       \setboxz@h{\refsfont@\keyformat{#2}}%
       \refindentwd\wd\z@
   \else\ifx#1\page\setbox\z@\hbox{\quad\bf#2}%
       \pagenumwd\wd\z@
   \else\ifx#1\item
       \setboxz@h{#2}\rosteritemwd=\wd\z@
   \else\ifx#1\itemitem
       \setboxz@h{#2}\rosteritemitemwd=\wd\z@
	\advance\rosteritemitemwd by .5pc
	\advance\rosteritemitemwd by \rosteritemwd
   \else\ifx#1\itemitemitem
       \setboxz@h{#2}\rosteritemitemitemwd=\wd\z@
	\advance\rosteritemitemitemwd by .5pc
	\advance\rosteritemitemitemwd by \rosteritemitemwd
   \else\message{\string\widestnumber\space not defined for this
      option (\string#1)}%
\fi\fi\fi\fi\fi\fi\fi\fi\fi}

\refstyle{A}
\widestnumber\key{M}

\endgroup

\catcode`\@=13


\catcode`\@=11

\def\journame{Contemporary Mathematics}

\issn{0271-4132}

\catcode`\@=13


\newskip\vadjustskip \vadjustskip=-\normalbaselineskip
\def\centertext
 {\hoffset=\pgwidth \advance\hoffset-\hsize
  \advance\hoffset-2truein \divide\hoffset by 2\relax
  \voffset=\pgheight \advance\voffset-\vsize
  \advance\voffset-2truein \divide\voffset by 2\relax
  \advance\voffset\vadjustskip
 }
\newdimen\pgwidth\newdimen\pgheight
\def\letter{letter}\def\AFour{AFour}
\ifx\PaperSize\letter
 \pgwidth=8.5truein \pgheight=11truein
 \message{- Got a paper size of letter.  }\centertext
\fi
\ifx\PaperSize\AFour
 \pgwidth=210truemm \pgheight=297truemm
 \message{- Got a paper size of AFour.  }\centertext
\fi
\def\GetNext#1 {\def\NextOne{#1}\if\relax\NextOne\let\next=\relax
        \else\let\next=\DoIt \fi \next}
\def\DoIt{\Act\NextOne\GetNext}
\def\ActOn#1{\expandafter\GetNext #1\relax\ }
\def\defcs#1{\expandafter\xdef\csname#1\endcsname}

 \font\twelvebf=cmbx12          
 \font\smc=cmcsc10              
\catcode`\@=11          
\def\eightpoint{\eightpointfonts
 \setbox\strutbox\hbox{\vrule height7\p@ depth2\p@ width\z@}%
 \eightpointparameters\eightpointfamilies
 \normalbaselines\rm
 }
\def\eightpointparameters{%
 \normalbaselineskip9\p@
 \abovedisplayskip9\p@ plus2.4\p@ minus6.2\p@
 \belowdisplayskip9\p@ plus2.4\p@ minus6.2\p@
 \abovedisplayshortskip\z@ plus2.4\p@
 \belowdisplayshortskip5.6\p@ plus2.4\p@ minus3.2\p@
 }
\newfam\smcfam
\def\eightpointfonts{%
 \font\eightrm=cmr8 \font\sixrm=cmr6
 \font\eightbf=cmbx8 \font\sixbf=cmbx6
 \font\eightit=cmti8
 \font\eightsmc=cmcsc8
 \font\eighti=cmmi8 \font\sixi=cmmi6
 \font\eightsy=cmsy8 \font\sixsy=cmsy6
 \font\eightsl=cmsl8 \font\eighttt=cmtt8}
\def\eightpointfamilies{%
 \textfont\z@\eightrm \scriptfont\z@\sixrm  \scriptscriptfont\z@\fiverm
 \textfont\@ne\eighti \scriptfont\@ne\sixi  \scriptscriptfont\@ne\fivei
 \textfont\tw@\eightsy \scriptfont\tw@\sixsy \scriptscriptfont\tw@\fivesy
 \textfont\thr@@\tenex \scriptfont\thr@@\tenex\scriptscriptfont\thr@@\tenex
 \textfont\itfam\eightit        \def\it{\fam\itfam\eightit}%
 \textfont\slfam\eightsl        \def\sl{\fam\slfam\eightsl}%
 \textfont\ttfam\eighttt        \def\tt{\fam\ttfam\eighttt}%
 \textfont\smcfam\eightsmc      \def\smc{\fam\smcfam\eightsmc}%
 \textfont\bffam\eightbf \scriptfont\bffam\sixbf
   \scriptscriptfont\bffam\fivebf       \def\bf{\fam\bffam\eightbf}%
 \def\rm{\fam0\eightrm}%
 }

\def\vfootnote#1{\insert\footins\bgroup
 \eightpoint\catcode`\^^M=5\leftskip=0pt\rightskip=\leftskip
 \interlinepenalty\interfootnotelinepenalty
  \splittopskip\ht\strutbox 
  \splitmaxdepth\dp\strutbox \floatingpenalty\@MM
  \leftskip\z@skip \rightskip\z@skip \spaceskip\z@skip \xspaceskip\z@skip
  \textindent{#1}\footstrut\futurelet\next\fo@t}

\def\p.{p.\penalty\@M \thinspace}
\def\pp.{pp.\penalty\@M \thinspace}
\newcount\sctno
\def\sctn#1\par
  {\removelastskip\vskip0pt plus4\normalbaselineskip \penalty-250
  \vskip0pt plus-4\normalbaselineskip \bigskip\medskip
  \centerline{\bf#1}\nobreak\medskip
}

\def\sct#1 {\sctno=#1\relax\sctn#1. }

\def\item#1 {\par\indent\indent\indent
 \hangindent3\parindent
 \llap{\rm (#1)\enspace}\ignorespaces}
 \def\inpart#1 {{\rm (#1)\enspace}\ignorespaces}
 \def\part {\par\inpart}

\def\Cs#1){\(\number\sctno.#1)}
\def\part#1 {\par\(#1)\enspace\ignorespaces}

\def\dsc#1 #2.{\medbreak{\bf\Cs#1)} {\it #2.} \ignorespaces}
\def\proclaim#1 #2 {\medbreak
  {\bf#1 (\number\sctno.#2)}\enspace \bgroup
\it}
\def\endproclaim{\par\egroup\medskip}

\def\dfn#1 {\medbreak {\bf Definition (\number\sctno.#1)}\enspace}
\def\rmk#1 {\medbreak {\bf Remark (\number\sctno.#1)}\enspace}
 \newcount\refno \refno=0        \def\NoKey{*!*}
 \def\MakeKey{\advance\refno by 1 \expandafter\xdef
  \csname\TheKey\endcsname{{\number\refno}}\NextKey}
 \def\NextKey#1 {\def\TheKey{#1}\ifx\TheKey\NoKey\let\next\relax
  \else\let\next\MakeKey \fi \next}
 \def\RefKeys #1\endRefKeys{\expandafter\NextKey #1 *!* }
 \def\SetRef#1 #2,{\hang\llap
  {\csname#1\endcsname.\enspace}{\smc #2},}
 \newbox\keybox \setbox\keybox=\hbox{25.\enspace}
 \newdimen\keyindent \keyindent=\wd\keybox
\def\references{\kern-\medskipamount
  \sctn \smc References\par
  \bgroup   \frenchspacing   \eightpoint
   \parindent=\keyindent  \parskip=\smallskipamount
   \everypar={\SetRef}\par}
\def\endreferences{\egroup}

 \def\serial#1#2{\expandafter\def\csname#1\endcsname ##1 ##2 ##3
        {\unskip\ {\it #2\/} {\bf##1} (##2), ##3}} 

\def\UThin{\penalty\@M \thinspace\ignorespaces}
\def\(#1){{\let~=\UThin\rm(#1)}}
\def\relaxnext@{\let\next\relax}
\def\cite#1{\relaxnext@
 \def\nextiii@##1,##2\end@{\unskip\space{\rm[\SetKey{##1},\let~=\UThin##2]}}%
 \in@,{#1}\ifin@\def\next{\nextiii@#1\end@}\else
 \def\next{{\rm[\SetKey{#1}]}}\fi\next}
\newif\ifin@
\def\in@#1#2{\def\in@@##1#1##2##3\in@@
 {\ifx\in@##2\in@false\else\in@true\fi}%
 \in@@#2#1\in@\in@@}
\def\SetKey#1{{\bf\csname#1\endcsname}}

\catcode`\@=12  

\let\:=\colon \let\ox=\otimes \let\x=\times
  \let\?=\overline

\def\smashedlongrightarrow{\setbox0=\hbox{$\longrightarrow$}\ht0=1pt\box0}
\def\risom{\buildrel\sim\over{\smashedlongrightarrow}}
 \def\lgto{-\mathrel{\mkern-10mu}\to}
 \def\smashedlgto{\setbox0=\hbox{$\scriptstyle\lgto$}\ht0=1.85pt
        \lower1.25pt\box0}

\def\tsum{\textstyle\sum}

  \def\es{^{\fam0 es}}  

\def\Act#1{\defcs{c#1}{{\fam2#1}}}               
 \ActOn{C I J K L M N O P Q }
 \def\Act#1{\defcs{#1}{\mathop{\fam0#1}\nolimits}} 
 \ActOn{cod sts frs frs$_1$ Hilb rts type wt Supp }
\def\Act#1{\defcs{c#1}{\mathop{\it#1}\nolimits}}
 \ActOn{Cok Ext Hom Ker Sym }
\def\Act#1{\defcs{#1}{\hbox{\rm\ #1 }}}         
 \ActOn{and by for where with on }
\def\Act#1{\defcs{I#1}{{\fam6#1}}}                
 \ActOn{A B C D G H P R S T V Z m }

 \RefKeys
 A76 At74 BL97 BL98 BL98s CH96 C90 C97 D73 E73 EC15 Gtt98 Gtz78 Grn88
GLS98 Hi53 Ho56 KPB K--P L98 NV97 Ran89 Sc68 Sh91 Sh97 S76 T76 V95 V97
W74 Z65
 \endRefKeys
\topmatter
\title Enumerating singular curves on surfaces\endtitle
\author Steven Kleiman and Ragni Piene\endauthor
\leftheadtext{S Kleiman and R Piene}%
\rightheadtext{Enumerating singular curves on surfaces}%

\address Mathematics Department, Room {\sl 2-278} MIT,
 {\sl77} Mass Ave, Cambridge, MA {\sl02139-4307}, USA\endaddress
\email kleiman\@math.mit.edu\endemail

\thanks The first author was supported in part by the US NSF and by the
Danish NSRC.\endthanks

\address Department of Mathematics, University of Oslo,
 PO Box {\sl1053} Blindern,  N-{\sl0316} Oslo, Norway\endaddress

 \email ragnip@math.uio.no\endemail

\issueinfo{241}
{}
{}
{1999}


\abstract
 We enumerate the singular algebraic curves in a complete linear system
on a smooth projective surface.  The system must be suitably ample in a
rather precise sense.  The curves may have up to eight nodes, or a
triple point of a given type and up to three nodes.  The curves must
also pass through appropriately many general points.  The number of
curves is given by a universal polynomial in four basic Chern numbers.

 To justify the enumeration, we make a rudimentary classification of the
types of singularities using Enriques diagrams, obtaining results like
Arnold's.  We show that the curves in question do appear with
multiplicity $1$ using the versal deformation space, Shustin's
codimension formula, and Gotzmann's regularity theorem.  Finally, we
relate our work to Vainsencher's work with up to seven nodes.\endabstract

\endtopmatter
\document

\sct1 Introduction

Consider a complete linear system of dimension $n$ on a smooth
irreducible complex projective surface.  Among the curves in the system,
some have $r$ ordinary nodes and no other singularities; they are, for
short, {\it r-nodal}.  How many $r$-nodal curves pass through $n-r$
general points?

Remarkably, at least for $r\le8$, if the system is suitably ample, then
the number $N_r$ of $r$-nodal curves is given by a {\it polynomial}
$P_r$ in the four basic Chern numbers of the surface $S$ and the system
$|\cL|$.  These Chern numbers are
 $$d:=\cL\cdot\cL,\ k:=\cL\cdot \cK_S,\ s:=\cK_S\cdot \cK_S,\ x:=c_2(S).$$
 A precise assertion is given in Theorem~(1.1) below.  It can be
extended to count curves with other types of singularities; one such
extension is Theorem~(1.2).  The proofs of these theorems is the subject
of this paper.

For example, $N_1$ is the number of nodal curves in a general pencil,
and it is given by the ``Zeuthen--Segre formula,"
        $$N_1=3d+2k+x.$$
 This formula was obtained by Hirzebruch \cite{Hi53} as a corollary of a
general formula for the index sum of a meromorphic vector field.  In our
case, the index sum is equal to $N_1+\cL\cdot\cL$, since $\cL\cdot\cL$
is just the number of base points.

In general, the polynomials $P_r$ may be obtained as follows.  Consider the
formal identity in $t$,
 $$\textstyle
 \sum_{r\ge0}P_rt^r/r!
        = \exp\bigl(\sum_{q\ge1}a_q t^q/q!\bigr).$$
 Thus $P_0=1$, and $P_1=a_1$, and $P_2= a_1^2+a_2$, and
$P_3=a_1^3+3a_2a_1+a_3$, and so forth. Set
 $$\eightpoint\eqalign{
a_1 &{}:= 3d+2k+x,\cr
a_2 &{}:=-42d-39k-6s-7x,\cr
a_3 &{}:= 1380d+1576k+376s+138x,\cr
a_4 &{}:= -72360d-95670k-28842s-3888x,\cr
a_5 &{}:= 5225472d+7725168k+2723400s+84384x,\cr
a_6 &{}:= -481239360d-778065120k-308078520s+7918560x,\cr
a_7 &{}:= 53917151040d+93895251840k+40747613760s-2465471520x,\cr
a_8 &{}:=-7118400139200d-13206119880240k-6179605765200s
        +516524964480x.\cr}$$
 View $P_r$ as a polynomial in $d,k,s,x$.  Then we have the following
theorem.
 \medskip
 {\bf Theorem (1.1)}\enspace {\it Assume $\cL=\cM^{\ox m}\ox \cN$
where $\cN$ is spanned and $\cM$ is very ample.  If $r\le8$ and $m\ge
3r$, then}
        $$N_r=P_r(d,k,s,x)/r!\,.$$
 \medskip
 Theorem (1.1) is extended in Theorem (1.2) below to enumerate curves
with a triple point of a given type and up to three nodes.
Specifically, let
        $$N(3),\ N(3,2),\ N(3,2,2),\and N(3,2,2,2)$$
 be the numbers of curves having an ordinary triple point and zero to
three ordinary nodes, and passing through $n-r$ general points where $r$
is 4, 5, 6, and 7 respectively.  In addition, let
        $$N(3(2)),\ N(3(2),2),\and N(3(2)')$$
 be the numbers of curves, respectively, (1) having a $D_6$ singularity,
that is, a triple point composed of a tacnode cut transversely by a
third branch, (2) having such a triple point and a distant ordinary
node, and (3) having an $E_7$ singularity, that is, a triple point
composed of a cusp cut tangentially by a second branch.  These curves
are also required to pass through $n-r$ general points where $r$ is 6,
7, and again 7, respectively. \medskip

  {\bf Theorem (1.2)} \enspace {\it As in Theorem \(1.1), assume
$\cL=\cM^{\ox m}\ox \cN$ where $\cN$ is spanned, $\cM$ is very ample,
and $m\ge 3r$ where $r$ is the appropriate number, now between $4$ and
$7$ as specified just above.  Then the following seven formulas hold:}
$$\eightpoint\def\&{\cr&\kern-2.5em}\eqalign{
     N(3)&{}=15d+20k+5s+5x;\cr
 N(3,2)&{}=45d^2+(15s+90k+30x-420)d+40k^2+(10s+30x-624)k\&
        +(5x-196)s+5x^2-100x;\cr
 N(3,2,2)&{}=\bigl(135d^3+(135x+45s+360k-3150)d^2+(300k^2+(60s+240x-6849)k\&
  +(-1476+30x)s+45x^2-1755x+18480)d+80k^3+(100x-3276+20s)k^2\&
  +((-1099+20x)s+40x^2-1983x+29946)k
  -30s^2+(10932-457x+5x^2)s\&+5x^3-235x^2+3120x\bigr)/2!;\cr
 }$$$$\eightpoint\def\&{\cr&\kern-2.5em}\eqalign{
 N(3,2,2,2)&{}=\bigl(405d^4+(-17010+1350k+135s+540x)d^3+(1620k^2
  +(270s+1350x\&-48573)k +(135x-7992)s+270x^2-14985x+239940)d^2
  +(840k^3+(180s+1080x\&-43074)k^2+((-11691+180x)s+450x^2-29052x+559398)k
  -270s^2+(-5013x\&+143184+45x^2)s+60x^3-4320x^2+113910x-1135080)d
  +160k^4+(40s+280x\&-12168)k^3+((-4242+60x)s +180x^2-13038x+284204)k^2
  +(-180s^2+(115156\&+30x^2-3687x)s+50x^3-4287x^2+144002x-1977552)k
  +(-90x+5408)s^2\&+(5x^3-783x^2+41282x-807006)s
  +5x^4-405x^3+12150x^2-128700x\bigr)/3!;\cr
 N(3(2)) &{}= 28x+168s+224d+406k;\cr
 N(3(2),2) &{}= -546x-7281s-8316d-16008k+28x^2+462xk+308xd\&
  +168sx+336sk+504sd+812k^2+1666kd+672d^2;\cr
 N(3(2)') &{}= 252d+488k+217s+42x.\cr}$$
 \medbreak

These seven formulas were chosen for a reason: they are explicitly used
as correction terms in Vainsencher's determination of $N_r$ for $r\le7$,
although they are not explicitly given in his paper \cite{V95}.
Vainsencher's treatment inspired ours greatly, but ours is more refined:
our formulas are more compact; our notion of ``suitably ample'' is more
precise; and our supplementary condition is to pass through $n-r$
general points, not simply to lie in a ``suitable general'' subsystem of
dimension $n-r$.  We also go further: Vainsencher was unsure about
$N_7$; we settle $N_7$ and treat $N_8$.  The seven formulas above are
only implicit in our proof of Theorem (1.1), but in Section~4 we explain
how to make them explicit to prove Theorem (1.2).  In Section~5, we compare
and contrast in more detail the nature of the appearance of the formulas
in Vainsencher's work with that in ours.

It is possible, by modifying our treatment, to enumerate the curves
having any other equisingularity type and passing through the
appropriate number $n-r$ general points, at least if $r\le8$.  Although
we don't do this enumeration here, nevertheless we must classify the
types of singularities involved to establish the validity of the
formulas in Theorems (1.1) and (1.2).  In Section~2, we carry out this
classification using Enriques diagrams; the results agree with Arnold's
classification \cite{A76} (see also Siersma's paper \cite{S76}).

To establish the validity of the formulas, we must consider the curves
lying in the system and having a given equisingularity type: we must
show that these curves form a reduced subsystem of the appropriate
codimension.  G\"ottsche \cite{Gtt98, 5.2} took a first step in this
direction, treating nodes in an ad hoc fashion, and his work helped
inspire ours.  Notably, it is his idea to take $\cL$ to be of the form
$\cM^{\ox m}\ox \cN$, and not simply of the form $\cM^{\ox m}$.
Moreover, he \cite{Gtt98, 2.2} suggested that it might be sufficient to
take $m\ge Cr$ for some universal constant $C$; in fact, our theorems
assert that it is sufficient to take $C=3$.  In Section~3, we study
equisingular systems of curves by using the theory of complete ideals
and the theory of the versal deformation space, supplemented by our
classification in Section~2 and by Shustin's codimension formula
\cite{Sh91, Thm.}.

We treat the eight formulas of Theorem (1.1), but not the seven formulas
of Theorem (1.2), in \cite{K--P}.  In Section 4 below, we review that
treatment for two reasons.  First, we can then explain how to modify it
to obtain the seven formulas of Theorem (1.2).  Second, we can then
explain how to use the results we prove in Sections 2 and 3 to obtain an
alternative proof of the validity of all the formulas.  In \cite{K--P},
with an eye toward any $r$, we establish validity by proving some
general results about Enriques diagrams, and by making a direct study of
the Hilbert scheme of clusters of points.  Whereas here, we apply
Gotzmann's regularity theorem \cite{Gtz78}; there, we prove an ad hoc
regularity result.  Here, we work only in characteristic zero; there, we
establish the validity in any characteristic for $r\le8$ if
$\cL=\cM^{\ox m}\ox \cN$ where $\cM$ and $\cN$ are spanned, $\cM$ is
only ample, and $m\ge m_0$ where
     $$m_0:= 3r+g^2+g+4-(s+x)/12\and g:=1+\cM\cdot(\cK_S+\cM)/2.$$
 However, in positive characteristic $p$, in a given count on an
irreducible surface, all the curves may possibly appear with the same
multiplicity $p^e$ where $e\ge1$.

There appears to be, unfortunately, no compact description of the sum,
        $$N(3)+N(3,2)t+N(3,2,2)t^2+N(3,2,2,2)t^3,$$
 like the compact exponential description given above for
$\sum_rN_rt^r$.  (The polynomial giving $N_6$ fills half a page in
\cite{V95, p.~514}; that giving $N_8$ fills a page and a half when
expanded.)  However, the displayed sum is the degree of a cycle class on
a suitable space, and this class has an analogous exponential
description.  The latter is explained in Section~4.  We give some other
examples of this phenomenon in \cite{KPB}.  We conjecture that
this cycle-theoretic exponential description holds for any suitably
general {\it algebraic\/} system of curves on any  {\it algebraic\/}
family of surfaces.

The exponential description of $\sum_rN_rt^r$ is not found in
Vainsencher's paper \cite{V95}.  Rather, it was discovered later by
G\"ottsche \cite{Gtt98, 2.3}.  He proved that, if each $N_r$ is given by
some polynomial for $r< r_0$ for some $r_0$, then necessarily
$\sum_{r<r_0} N_rt^r$ is of the form $\exp\bigl(\sum_{q< r_0}
A_qt^q/q!\bigr)$ where each $A_q$ is some {\it linear\/} combination of
$d,k,s,x$.  This linearity reflects another remarkable property of the
polynomials: they continue to work when $S$ is replaced by a surface
with several (disjoint) components.  (In the summer of 1997 when
G\"ottsche wrote \cite{Gtt98}, considerable evidence suggested that
$r_0=\infty$.  Nine months later, Liu \cite{L98} offered a symplectic
proof, and in February 1999, he claimed [pvt. comm.] to have generalized
his result to any equisingularity type.)  Our use of the exponential was
inspired by G\"ottsche's discovery, but is logically independent of it.

G\"ottsche related his $A_q$ to the Fourier developments of certain
modular forms.  These relations suffice to determine the $A_q$
completely when $\cK_S$ vanishes.  For a K3 surface and for an Abelian
surface, the resulting formulas agree with the formulas proved by Bryan
and Leung \cite{BL97}, \cite{BL98}, \cite{BL98s} (see also the authors'
paper \cite{KPB}).  Furthermore, G\"ottsche explained how to compute the
$A_q$ from recursive formulas, such as those of Caporaso and Harris
\cite{CH96} for $\IP^2$, of Ran \cite{Ran89} for $\IP^2$, or of Vakil
\cite{V97} for $\IP^1\x\IP^1$.  For $q\le8$, he checked that the
resulting $A_q$ are equal to the $a_q$ in Theorem (1.1) above, which
were obtained independently (they were also obtained from scratch,
except for one unknown in $a_8$).

In short, in Section 2, we make an abstract combinatorial classification
of the Enriques diagrams associated to the curves of interest to us, and
tabulate the values of their basic numerical characters.  In Section 3,
we interpret these characters geometrically, and analyze the locus of
curves with a given diagram, or equivalently, a given equisingularity
type.  In Section 4, we derive the formulas and establish their
validity; this material rests in part on \cite{K--P}.  Finally, in
Section 5, we relate Vainsencher's treatment to ours.

\input xy \xyoption{all}
\sct2 Minimal Enriques diagrams

\def\lab#1#2#3{\relax
 \POS*=0\dir{*}="a",*+!#1=0{\labelstyle #2},"a",#3}
\def\rs#1#2{\relax\kern-1.75em\raise1.75pc\hbox{$#1_#2$}}
\def\us#1#2 {\relax\kern-1.75em\raise1.75pc\hbox{$#1_{#2}$}}
\def\fig#1 #2\par{\par\kern6pt{\advance\leftskip by 1.5\parindent
 \noindent{\bf Fig.\thinspace 2-#1. }Diagram #2\par\medbreak}}

In this section we are going to make a rudimentary combinatorial
classification of the ``minimal Enriques diagrams'' that arise from the
curves of special interest to us.  The diagrams represent the
equisingularity types of the curves (this fact was discovered by
Enriques \cite{EC15, IV.I}, and proved rigorously by Zariski \cite{Z65},
although Zariski did not use the language of diagrams; however, this
language was used by Casas \cite{C97, \p.100}).

Given a reduced curve on a smooth surface (which need not be complete),
form the configuration of all infinitely near points on all the branches
of the curve through all its singular points.  Weight each infinitely
near point with its multiplicity on the strict transform of the curve.
Also, say that one infinitely near point is ``proximate'' to a second if
the first lies above the second and on its strict transform, that is, on
the strict transform of the exceptional divisor of the blowup centered
at the second.  Equipped with this binary relation of proximity, the
weighted configuration has an abstract combinatorial structure, which
Enriques \cite{EC15, IV.I, \pp.350--51} encoded in a convenient diagram.

By the theorem of strong embedded resolution, all but finitely many
infinitely near points are of multiplicity 1, and are proximate only to
their immediate predecessors.  So it is common to prune off all the
infinite unbroken successions of such points, leaving finitely many
points.  They are known as the {\it essential} points of the curve (in
the terminology of \cite{GLS98, 2.2}).  Thus, the essential points form
a configuration, whose end points either are of higher multiplicity or
are proximate to a remote predecessor.  We will call the abstract
combinatorial structure of this configuration the ``minimal Enriques
diagram'' of the curve.  (In \cite{C97, 3.9}, Casas adds all the
successors of the remote end points, although these successors are free
points of weight 1, and he calls the resulting combinatorial structure
the ``Enriques diagram'' of the curve.)

We formally define a {\it minimal Enriques diagram\/} $\ID$ to be a
finite weighted forest that is equipped with a binary relation and that
is subject to the laws stated below.  As usual, a {\it forest\/} is a
disjoint union of trees.  A {\it tree\/} is a directed graph, without
loops; it has a single initial vertex, or {\it root,} and every other
vertex has a unique immediate predecessor.  A final vertex is called a
{\it leaf}.

For convenience, we consider a vertex to be one of its own predecessors
and one of its own successors, and we call its other predecessors and
successors {\it proper}.  We require the {\it weight\/} of a vertex $V$
to be an integer $m_V$ at least 1.  We denote the binary relation by
`$\succ$' and call it {\it proximity}.  As is customary, if a vertex $V$
is proximate to a remote predecessor, then we call $V$ a {\it
satellite\/}; otherwise, we call $V$ {\it free}.  Thus a root is free.

Furthermore, $\ID$ is subject to the following laws:
 \smallskip{\advance\leftskip by 2\parindent\it\parindent=0pt
 \(Law of Proximity)\enspace A root is proximate to no vertex.  If
a vertex is not a root, then it is proximate to its immediate
predecessor and to at most one other vertex; the latter must be a remote
predecessor.  If one vertex is proximate to a second, and if a third
lies properly between the two, then it too is proximate to the second.
 \smallskip
 \(Proximity Inequality)\enspace  For each vertex $V$,
        $$m_V\ge\tsum_{W\succ V}m_W.$$
 \smallskip
 \(Law of Succession)\enspace A vertex $V$ may have at most $m_V$
immediate successors, of which any number may be free, but at most two
may be satellites, and they must be satellites of different vertices.
 \smallskip
 \(Law of Minimality)\enspace Every leaf of weight $1$ is a satellite.
 \smallskip}

After Enriques, we depict $\ID$ graphically as follows.  We shape the
sequence of edges connecting a maximal succession of free vertices into
a smooth curve.  We shape the sequence of edges connecting a maximal
succession of vertices that are all proximate to the same vertex, $T$
say, into a line segment.  Its first vertex, $V$ say, is an immediate
successor of $T$.  The edge from $T$ to $V$ joins the segment at an
angle.  Specific diagrams are given in the figures below.

For convenience, set
  $$\eqalign{
         \frs(\ID)&:=\hbox{the number of free vertices in $\ID$,}\cr
         \rts(\ID)&:=\hbox{the number of roots in $\ID$.}\cr}$$
 In terms of the preceding numbers, we define the following characters:
  $$\eqalign{
        \dim(\ID)&:=\rts(\ID)+\frs(\ID);\cr
        \deg(\ID)&:=\tsum_{V\in\ID}{m_V+1\choose2};\cr
        \cod(\ID)&:=\deg(\ID)-\dim(\ID);\cr}\quad
 \eqalign{\delta(\ID)&:=\tsum_{V\in\ID}{m_V\choose2};\cr
        r(\ID)&:=\tsum_V \bigl(m_V-\tsum_{W\succ V}m_W\bigr);\cr
        \mu(\ID)&:=2\delta(\ID)-r(\ID)+\rts(\ID).\cr}$$
 Each of these characters has a geometric meaning, which will be
explained in the next section.  Note in passing that $r(\ID)$ measures
the total failure of the proximity inequalities to be equalities.

In the next section, we'll need to know which minimal Enriques diagrams
$\ID$ have $\cod(\ID)\le9$.  It suffices to know the $\ID$ having only
one root $R$ since the others are disjoint unions of these, and so we'll
classify them.  We'll find that $m_R$ must be 2, 3, or 4, and while
we're at it, we'll identify all the $\ID$ where $m_R$ is 2 or 3 and
$\cod(\ID)$ is arbitrary, and all the $\ID$ where $m_R$ is 4 and
$\cod(\ID)\le10$.  Our results are given in the figures below.

The figures are labelled according to Arnold's scheme \cite{A76}.  He
classified a large number of planar singularities, or bivariate analytic
function germs, and gave normal forms for each topological equivalence
class.  His results imply ours (in characteristic zero), but it is
rather easy and elementary to obtain ours directly, as we do now.

Fix a minimal Enriques diagram $\ID$ with one root $R$.  We begin the
classification by proving a useful little lemma.
 \medbreak
 {\bf Lemma (2.1)}\enspace {\it Let $S$ be a vertex.  Let $V$ and $W$
be successors of $S$.  Then $m_R\ge m_S\ge m_V$.  If $m_S= m_V$, then $W$
is either a successor of\/ $V$ or a predecessor of\/ $V$; if the latter,
then $m_W= m_V$, and $V$ is not remotely proximate to $W$.  If $m_S=1$
and $V$ is a leaf, then $m_V=1$ and $V$ is remotely proximate to some
proper predecessor of $S$.}
 \medskip
  Indeed, if $V=R$, then the first two assertions are trivial, and the
third is vacuous because its hypotheses imply that $R$ is a leaf of
weight 1 in violation of the Law of Minimality.  So assume that $V\neq
R$, and let $U$ be the immediate predecessor of $V$.  Then by the Law
of Proximity and the Proximity Inequality, $m_U\ge m_V$, and if $m_U=
m_V$, then $V$ is the only vertex proximate to $U$.  Since $U$ is closer
to $R$ than $V$ is, we may assume by induction that the lemma holds with
$U$ as $V$.  Hence, $m_R\ge m_U\ge m_V$.  If $S=V$, then the first
assertion follows, the second is trivial, and the third holds by the Law
of Minimality.  So assume that $S\neq V$.  Then $S$ is a predecessor of
$U$.  Hence, $m_R\ge m_S\ge m_U$ by induction.  Therefore, $m_R\ge
m_S\ge m_V$.

Suppose $m_S= m_V$.  Then $m_S= m_U=m_V$.  Hence, by induction, $W$ is
either a successor of $U$ or a predecessor of $U$; if the latter, then
$m_W= m_U$.  If $W$ is a proper successor of $U$, then it must be a
successor of $V$; otherwise, $U$ would have two immediate successors,
and both would be proximate to $U$ by the Law of Proximity, contrary to
what we observed above.  So assume that $W$ is a predecessor of $U$.
Then $m_W= m_U=m_V$.  Moreover, $V$ is not remotely proximate to $W$.
Otherwise, $W\neq U$ since $U$ is the immediate predecessor of $V$.
Hence $U$ too would be proximate to $V$ by the Law of Proximity.
Therefore, the Proximity Inequality would yield
        $$m_W\ge m_U+m_V=2m_W.$$

Finally, assume $m_S=1$ and $V$ is a leaf.  Then $m_V=1$ since $m_S\ge
m_V$.  By the Law of Minimality, $V$ is proximate to some remote
predecessor, and it cannot be a successor of $S$ by the second
assertion.  The lemma is now proved.

The last assertion of the lemma implies that $m_R\neq 1$.

Assume $m_R=2$.  Let $L$ be a leaf.  By the lemma, $2=m_R\ge m_L$.  If
$m_L=2$, then the lemma implies that $\ID$ is the diagram $A_{2i-1}$
described in Fig.~2-1.
 \vskip-\smallskipamount
$$\xymatrix @*=0{%
&&&\lab{DL}{2}{}&&&&\lab{DL}{2}{}\\
        \rs A1 \lab{UL}{2}{}
        &&\rs A3\lab{UL}{2}{\ar@{-}@/^/[ur]}
        &&&\rs A5\lab{UL}{2}{\ar@{-}@/^/[urr]^*\dir{*}_{2}}}
 $$
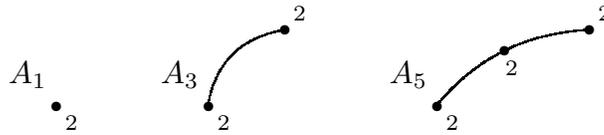
 \fig1 $A_{2i-1}$ for $i\ge1$: a succession of $i$ free vertices of
weight $2$

Suppose $m_L=1$.  Then, by the Law of Minimality, $L$ is proximate to
some remote predecessor $T$.  By the lemma, $2=m_R\ge m_T$.  Let $U$ be
the immediate predecessor of $L$.  Then, by the Law of Proximity, $U$
is also proximate to $T$.  So, by the Proximity Inequality, $m_T\ge
m_U+m_L$.  Hence $m_T=2$ and $m_U=1$.  Let $T'$ be the immediate
predecessor of $U$.  Then $T'=T$; otherwise, similarly,
        $$2=m_T\ge m_{T'}+m_U+m_L=3.$$
 The lemma now implies that $\ID$ is the diagram $A_{2i}$ described in
Fig.~2-2.
 \vskip-\smallskipamount
 $$\xymatrix @*=0{%
&\lab{DL}{1}{\ar@{-}[d]}&&&&\lab{DL}{1}{\ar@{-}[d]}\\
 \rs A2\lab{UL}{2}{\ar@{-}@/^/[ur]}
 &\lab{UL}{1}{}&&\rs A4\lab{UL}{2}{\ar@{-}@/^/[urr]^*\dir{*}_{2}}
 &&\lab{UL}{1}{}
 }{\qquad}$$
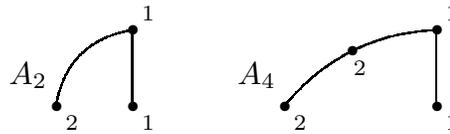
 \fig2 $A_{2i}$ for $i\ge1$: a succession of $i$ free vertices of
weight $2$, followed by two vertices of weight 1 both proximate to the
$i$th vertex

Assume $m_R=3$.  If $\ID=\{R\}$, then $\ID$ is the diagram $D_4$, shown
in Fig.~2-3.
 \vskip-\smallskipamount
 $$\xymatrix @*=0{%
&&&\lab{DL}{1}{\ar@{-}[d]}&&&\lab{DL}{2}{}&&&&\lab{DL}{1}{\ar@{-}[d]}\\
\rs D4\lab{UL}{3}{}&&\rs D5\lab{UL}{3}{\ar@{-}@/^/[ur]*{}}
 &\lab{UL}{1}{}&&\rs D6\lab{UL}{3}{\ar@{-}@/^/[ur]}
&&&\rs D7\lab{UL}{3}{\ar@{-}@/^/[urr]^*\dir{*}_{2}}&&\lab{UL}{1}{}
 }$$
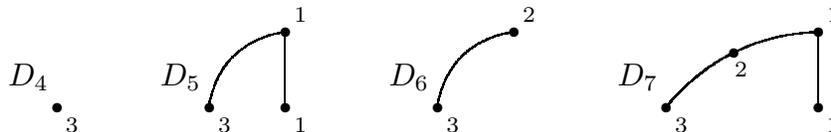
 \fig3 $D_k$ for $k\ge4$: a root of weight $3$ followed by nothing if
$k=4$, followed by two vertices of weight $1$ both proximate to the root
if $k=5$, and followed by the diagram $A_{k-5}$ if $k\ge6$

Suppose $R$ has an immediate successor $S$.  Then by the lemma $3=m_R\ge
m_S$.  First, suppose $m_S=1$.  Let $L$ be a leaf that succeeds $S$.
The lemma implies that $m_L=1$ and $L$ is remotely proximate to $R$.
Then $L\neq S$ because $S$ is an immediate predecessor of $R$.

Conceivably, $R$ has a second immediate successor $W$.  If so, then, by
the Law of Proximity and the Proximity Inequality, we'd have
        $$3=m_R\ge m_W+m_S+m_L=m_W+2\ge3.$$
 So $m_W=1$.  However, with $W$ for $S$, the argument above would imply
that $W$ has a proper successor $M$ that is remotely proximate to $R$.
So, by the Proximity Inequality, we'd have
        $$3=m_R\ge m_W+m_M+m_S+m_L=4.$$
 Thus $W$ does not exist.

If $L$ is an immediate successor of $S$ and the only one, then $\ID$ is
the diagram $D_5$, shown in Fig.~2-3.  Suppose not.  Then, by the lemma,
there is some vertex $V$ strictly between $S$ and $L$.  Then, by the
Law of Proximity, $V$ too is proximate to $R$.  If there were a second
such vertex $W$, then by the Proximity Inequality, we'd have
        $$3=m_R\ge m_S+m_V+m_W+m_L=4.$$
 Hence $\ID$ is the diagram $E_6$, shown in Fig.~2-4.
 \vskip-\smallskipamount
 $$\xymatrix @*=0{%
&\lab{DL}{1}{\ar@{-}[d]}&&&\lab{DL}{2}{\ar@{-}[d]}
&&&\lab{DL}{2}{\ar@{-}[d]}\\
\rs E6\lab{UL}{3}{\ar@{-}@/^/[ur]}&\lab{DL}{1}{\ar@{-}[d]}
&&\rs E7\lab{UL}{3}{\ar@{-}@/^/[ur]}&\lab{UL}{1}{}
&&\rs E8\lab{UL}{3}{\ar@{-}@/^/[ur]}
&\lab{UL}{1}{\ar@{-}[r]}&\lab{UL}{1}{}\\
&\lab{UL}{1}{}
}$$
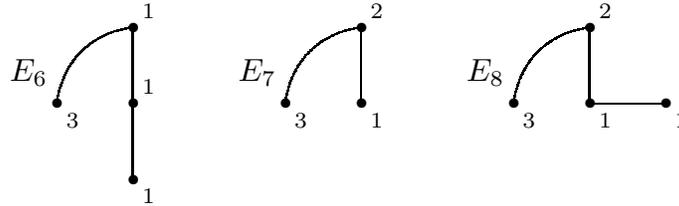
 \fig4 $E_{6l+j}$ for $l\ge1$ and $j=0,1,2$: a succession of $l$
vertices of weight $3$ followed by three vertices of weight $1$ all
proximate to the the $l$th vertex if $j=0$; followed by two vertices,
one of weight $2$ and one of weight $1$, both proximate to the the
$l$th vertex if $j=1$; and followed by three vertices, one of weight
$2$ and then two of weight $1$, the first two both proximate to the
$l$th vertex, and the second two both proximate to the $(l+1)$th
vertex if $j=2$

Second, suppose $m_S=2$.  Suppose there is no other vertex proximate to
$R$.  Then removing $R$ from $\ID$ produces a new minimal Enriques
diagram.  It has $S$ as its only root.  So it is the diagram $A_l$ for a
suitable $l$ by the case $m_R=2$ above.  Hence $\ID$ is the diagram
$D_{l+5}$ described in Fig.~2-3.

Suppose there is, other than $S$, another vertex $T$ proximate to $R$.
Then, by the Proximity Inequality,
        $$3=m_R\ge m_S+m_T=2+m_T\ge3.$$
 So $m_T=1$.  Arguing as before, we find that $T$ cannot be a second
immediate successor of $R$; otherwise, $T$ would have a proper successor
$M$ that is remotely proximate to $R$, and so, by the Proximity
Inequality,
        $$3=m_R\ge m_T+m_M+m_S=4.$$
 Therefore, $T$ is a successor of $S$.

Other than $S$, there is no immediate successor $T'$ of $R$, because
$T'$ would be proximate to $R$ by the Law of Proximity, but would not
be a successor of $S$, contrary to the preceding conclusion.

Other than $S$ and $T$, no vertex $T'$ is proximate to $R$; otherwise,
        $$3=m_R\ge m_S+m_T+m_{T'}\ge4.$$
 By the Law of Proximity, any vertex $T'$ between $S$ and $T$ is
proximate to $R$.  Hence $T$ is an immediate successor of $S$.  If $\ID$
has no further vertices, then it is the diagram $E_7$ shown in Fig.~2-4.

Conceivably, $S$ has a second immediate successor $T'$.  Arguing as
before, suppose so.  Then, by the Law of Proximity and the Proximity
Inequality,
        $$2= m_S\ge m_T+m_{T'}\ge2.$$
 So $m_{T'}=1$.  However, then $T'$ would have a proper successor $M'$
that is remotely proximate to $S$.  So, by the Proximity Inequality,
we'd have
        $$2=m_S\ge m_{T'}+m_{M'}+m_T=3.$$
 Thus $T'$ does not exist.

Suppose $\ID$ has a vertex $U$ in addition to $R$, $S$ and $T$.  Then,
by what we've just seen, $U$ must be a successor of $T$.  Let $L$ be a
leaf that succeeds $U$.  By the lemma, $L$ is remotely proximate to some
proper predecessor $S'$ of $U$.  Then $S'=S$ because $S$ and $T$ are the
only vertices proximate to $R$.  Therefore, $U=L$; otherwise, the
Proximity Inequality would yield
        $$2=m_S\ge m_T+m_U+m_L=3.$$
  Thus $\ID$ is the diagram $E_8$, shown in Fig.~2-4.

Third and finally, suppose $m_S=3$.  Then there is no other vertex
proximate to $R$ because of the Proximity Inequality.  Hence removing
$R$ from $\ID$ produces a new minimal Enriques diagram.  It has one less
vertex of weight 3, and it has $S$ as its only root.  So either it is
one of the diagrams already identified, or removing $S$ produces yet
another minimal Enriques diagram.  Continuing, we conclude that $\ID$ is
either the diagram $E_k$, described in Fig.~2-4, or the diagram
$J_{l,k}$, described in Fig.~2-5.
 \vskip-\smallskipamount
$$\xymatrix @*=0{%
&\lab{DL}{3}{}&&&&\lab{DL}{1}{\ar@{-}[d]}\\
 \us J2,0 \lab{UL}{3}{\ar@{-}@/^/[ur]}
 &&&\us J2,1 \lab{UL}{3}{\ar@{-}@/^/[urr]^*\dir{*}_{3}}&&\lab{UL}{1}{}
}$$
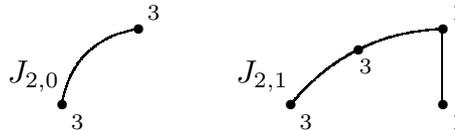
 \fig5 $J_{l,k}$ for $l\ge2$ and $k\ge0$: a succession of $l-1$
vertices of weight 3 followed by the diagram $D_{k+4}$

Lastly, assume $m_R\ge4$ and $\cod(\ID)\leq 10$.  In general, the
defining formula for $\cod(\ID)$ may be rewritten as
        $$\textstyle\cod(\ID)={m_R+1\choose2}-2
 +\tsum_{V\in\{\ID-R\}}\bigl({m_V+1\choose2}
        -\tau_V\bigr)$$
 where $\tau_V=1$ if $V$ is free, and $\tau_V=0$ if $V$ is a satellite.
Hence, $m_R=4$ since $\cod(\ID)\leq 10$.  Moreover, $\ID$ has no
additional free vertex of weight 3 or more.  Furthermore, either $\ID$
has no additional free vertex of weight 2 or more, and at most two
satellites of weight 1, or else $\ID$ has one free vertex of weight $2$
and no satellites.

A free vertex $V$ of weight 1 does not contribute to $\cod(\ID)$.
However, arguing much as before using the lemma, we find that $V$ is
followed by a nonempty succession of vertices of weight 1, all of which
are remotely proximate to a predecessor $T$ of $V$.  Necessarily, $T$ is
the immediate predecessor of $V$ because $V$ is free, yet proximate to
$T$ too.  It follows that $\ID$ is one of the diagrams shown in
Fig.~2-6, Fig.~2-7, and Fig.~2-8.

$$\xymatrix @*=0{%
&&&\lab{DL}{1}{\ar@{-}[d]}&&&\lab{DL}{2}{}\\
\us X1,0 \lab{UL}{4}{}&&\us X1,1 \lab{UL}{4}{\ar@{-}@/^/[ur]*{}}&\lab{UL}{1}{}%
&&\us X1,2 \lab{UL}{4}{\ar@{-}@/^/[ur]}
 }$$
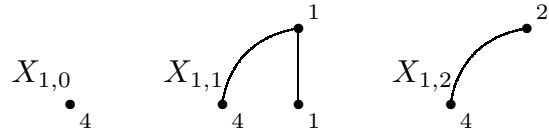
 \fig6 $X_{1,k}$ for $k=0,1,2$: a root of weight $4$ followed by
nothing if $k=0$, followed by two vertices of weight $1$ both proximate
to the root if $k=1$, and followed by one vertex of weight 2 if $k=2$

$$\xymatrix @*=0{%
\lab{DL}{1}{\ar@{-}[r]}&\lab{DL}{1}{}\\
&&\lab{DL}{1}{\ar@{-}[d]}\\
\us Y1,1 \lab{UL}{4}{\ar@{-}@/^/[urr],\ar@{-}@/_/[uur]}&&\lab{UL}{1}{}
 }$$
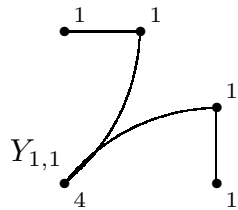
 \fig7 $Y_{1,1}$: a root of weight $4$ followed by two immediate
successors, each of weight $1$ and each having one successor of weight
$1$ that is proximate to the root

$$\xymatrix@*=0{%
&\lab{DL}{1}{\ar@{-}[d]}
\\
\us Z11 \lab{UL}{4}{\ar@{-}@/^/[ur]}&\lab{DL}{1}{\ar@{-}[d]}
\\
&\lab{UL}{1}{}
}$$
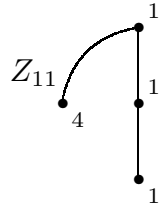
 \fig8 $Z_{11}$: a root of weight $4$ followed by three vertices of
weight $1$ all proximate to the root

Our classification is now complete.  A simple calculation yields all the
corresponding numerical characters, and they are listed in Table 2-1.
 Finally, here are all the minimal diagrams $\ID$ with $\cod(\ID)\leq
10$ and one root:
 $$A_1,\dots , A_{10}, D_4,\dots, D_{10}, J_{2,0}, J_{2,1}, E_6, E_7, E_8,
X_{1,0}, X_{1,1}, X_{1,2}, Z_{11}, Y_{1,1}.$$
 It is interesting to see the normal forms of the equations of curves
with these as associated diagrams.  So although we won't need them,
we've taken them from Arnold's paper  \cite{A76} and listed them in
Table 2-2.

Note in passing that we have not used the Law of Succession.  The first
time that it is needed is to rule out the following possible $\ID$ with
$\cod(\ID)=12$: a root of weight 4, followed by a vertex of weight 2,
followed by two immediate successors, each being of weight 1 and
proximate to the root.
 \vfill
$$\matrix \multispan7\hfil\bf Table \number\sctno-1\hfil\cr
\multispan7\hfil\bf Characters of Enriques diagrams\hfil\strut\cr
\noalign{\smallskip\hrule\smallskip}
\type & \cod & \deg & \dim & r &\delta & \mu \strut\cr
\noalign{\smallskip\hrule\smallskip}
 A_{2i-1} & 2i-1 & 3i & i+1 & 2 & i & 2i-1 \strut\cr
 A_{2i} & 2i & 3i+2 & i+2 & 1 & i & 2i \cr
 D_{2i} & 2i & 3i & i & 3 & i+1 & 2i \cr
 D_{2i+1} & 2i+1 & 3i+2 & i+1 & 2 & i+1 & 2i+1 \cr
 J_{l,2i}& 2i-1+5l & 3i+6l & i+1+l & 3 & i+3l & 2i-2+6l \cr
 J_{l,2i+1} & 2i+5l & 3i+2+6l & i+2+l & 2 & i+3l & 2i-1+6l \cr
 E_{6l} & 1+5l & 3+6l & l+2 & 1 & 3l & 6l \cr
 E_{6l+1} & 2+5l & 4+6l & l+2 & 2 & 1+3l & 1+6l \cr
 E_{6l+2} & 3+5l & 5+6l & l+2 & 1 & 1+3l & 2+6l \cr
 X_{1,0} & 8 & 10 & 2 & 4 & 6 & 9 \cr
 X_{1,1} & 9 & 12 & 3 & 3 & 6 & 10 \cr
 X_{1,2} & 10 & 13 & 3 & 4 & 7 & 11 \cr
 Z_{11} & 10 & 13 & 3 & 2 & 6 & 11 \cr
 Y_{1,1} &10 & 14 & 4 & 2 & 6 & 11 \cr
 \endmatrix$$
 \bigskip

$$\matrix\multispan4\hfil\bf Table 2-2\hfil\cr
\multispan4\hfil\bf Normal forms\hfil\strut\cr
\noalign{\smallskip\hrule\smallskip}
\type & \hbox{normal form}
 & \type & \hbox{normal form} \strut\cr
\noalign{\smallskip\hrule\smallskip}
 A_k & y^2+x^{k+1}\ (k\ge 1)
 & E_8 & x^3+y^5 \cr
 D_k & x^2y+y^{k-1}\ (k\ge 4)& X_{1,0}
 & x^4+ax^2y^2+y^4\ (a^2\neq 4)\cr
 J_{2,0}& x^3+ax^2y^2+y^6\ (4a^3+27\neq 0)
 & X_{1,1} & x^4+x^2y^2+ay^5\ (a\neq 0) \cr
 J_{2,1} & x^3+x^2y^2+ay^7\ (a\neq 0)
 & X_{1,2} & x^4+x^2y^2+ay^6\ (a\neq 0) \cr
 E_6 & x^3+y^4
 & Z_{11} & x^3y+y^5+axy^4 \cr
 E_7 & x^3+xy^3
 & Y_{1,1} & x^5+ax^2y^2+y^5\ (a\neq 0) \cr
 \endmatrix$$

\sct3 Geometric interpretation

In the last section, we introduced six numerical characters of an
abstract minimal Enriques diagram $\ID$.  In this section, we'll
interpret the six geometrically.  First, we'll show that three of them
represent standard characters of any reduced curve $C$ with $\ID$ as its
associated diagram.  Then we'll relate the other three to the set of all
such $C$ in a given complete linear system $|\cL|$.

Fix, once and for all, a nonempty abstract minimal Enriques diagram
$\ID$, a smooth surface $S$, and an invertible sheaf $\cL$ on $S$.
After proving Proposition (3.1), we'll assume that $S$ is irreducible
and projective so that the complete linear system $|\cL|$ is
parameterized by a projective space $Y$.

Let $C$ be any reduced curve on $S$.  One standard character of $C$ is
its total {\it Milnor number} $\mu(C)$, which can be defined by the
formula,
        $$\mu(C):=\dim H^0(\cO_S/\cJ_{C,S}),$$
 where $\cJ_{C,S}$ is the Jacobian ideal of $C$ on $S$; the latter is just
the first Fitting ideal of the sheaf of differentials $\Omega^1_C$
viewed as an $\cO_S$-module.

Define two characters of $C$ using the normalization map $\nu\:C'\to C$.
First, set
 $$r(C):=\hbox{\rm the number of points of $C'$ in }
        \nu^{-1}\Supp(\nu_*\cO_{C'}/O_C);$$
in other words, $r(C)$ is the total {\it number of branches\/} of $C$
through all its singular points.  Second, set
        $$\delta(C):=\dim H^0(\nu_*\cO_{C'}/O_C);$$
 we may call $\delta(C)$ {\it the genus discrepancy\/} of $C$.

These three characters of $C$ are, according to the following
proposition, equal to the corresponding characters of $\ID$ if $\ID$ is
the diagram we associate to $C$ in the way described at the beginning of
Section 2. \medbreak

{\bf Proposition (3.1)}\enspace {\it If $C$ has $\ID$ as its associated
diagram, then}
 $$r(\ID)=r(C) \and \delta(\ID)=\delta(C) \and \mu(\ID)=\mu(C).$$
 \medbreak

Indeed, consider $r(\ID)$ and $r(C)$.  Both vanish if $C$ is smooth, or
equivalently, if $\ID$ is empty.  To prove that $r(\ID)=r(C)$, proceed
by induction on the number of vertices in $\ID$.  Both $r(\ID)$ and
$r(C)$ are ``additive'' in the set of singular points of $C$.  So we may
assume that $C$ has only one singular point, $z$ say.  Let $R$ be the
corresponding root of $\ID$.

Blowup the surface at $z$.  Let $E$ be the exceptional divisor, $C'$ be
the strict transform of $C$, and $z_1,\dots,z_n$ be the points of $C'$
on $E$.  Then the number of branches of $C$ through $z$ is equal to the
sum over $i$ of the number of branches of $C'$ through $z_i$.  The
latter number is 1 if $z_i$ is a simple point of $C'$.  There are two
groups of these points: those whose branch is transverse to $E$, say
$z_1,\dots,z_l$, and those whose branch is not, say $z_{l+1},\dots,
z_k$.

To find $l$, note that $m_R$ is equal to the intersection number of $C'$
with $E$, so to the sum, over all the infinitely near points $w$
proximate to $z$, of the multiplicity at $w$ of the strict transform of
$C$.  Each $w$ corresponds to a vertex $W$ in $\ID$, proximate to $R$,
unless $w$ is one of $z_1,\dots,z_l$.  Hence, we have
        $$m_R-\tsum_{W\succ R}m_W=l.$$

To find $k-l$, let $z_i$ be one of $z_{l+1},\dots, z_k$.  Say $z_i$
corresponds to $S_i$ in $\ID$.  Consider all the vertices that succeed
$S_i$.  Lemma (3.1) implies that they form a single line of succession,
say from $V_{i1}:=S_i$ to $V_{ij_i}$; moreover, each $V_{ij}$ is
proximate to $R$ as well as to its immediate predecessor, which is
$V_{i,j-1}$ if $j\ge2$; also $m_{V_{ij}}=1$ for every $i,j$.  Hence, we
find
  $$\tsum_{i=l+1}^k\tsum_j^{j_i}\bigl(m_{V_{ij}}
                -\tsum_{W\succ V_{ij}}m_W\bigr)
        = \tsum_{i=l+1}^k\bigl(\tsum_j^{j_i-1}(1-1)+1\bigr)=k-l.$$

The remaining points $z_{k+1},\dots,z_n$ are the singular points of
$C'$.  So they correspond to the roots of the minimal Enriques diagram
$\ID'$ of $C'$.  By induction, the number of branches of $C'$ through
$z_{k+1},\dots,z_n$ is equal to $r(\ID')$.  Hence
        $$r(C)=l+(k-l)+r(\ID').$$
  Now, $\ID'$ is obtained from $\ID$ by deleting $R$ and also all the
$V_{ij}$.  Therefore, $r(C)=r(\ID)$.

The equation $\delta(\ID)=\delta(C)$ is a version of a famous old
result, which is often attributed to Noether and Enriques.  For a modern
proof, see Th\'eor\`eme (2.11)(ii) on \p.19 of \cite{D73}.

The equation $\mu(\ID)=\mu(C)$ follows immediately from the other two,
together with the celebrated Milnor--Jung formula,
        $$\mu(C)=2\delta(C)-r(C)+s(C),$$
 where $s(C)$ denotes the number of (distinct) singular points of $C$,
which is equal to the number of roots $\rts(\ID)$.  For a proof of the
formula, see Corollary 6.4.3 on \p.205 in \cite{C97}.  For more
references, see Lipman's review of J.-J. Risler's proof in Math Reviews
[MR 46-5334].  The proposition is now proved. \medbreak

{}From now on, assume that $S$ is irreducible and projective.  Denote by
$Y$ the projective space parameterizing the complete linear system
$|\cL|$, and by $Y(\ID)$ the subset of points representing curves $C$
with $\ID$ as associated diagram.  The following proposition interprets
the character $\cod(\ID)$ in terms of the geometry of $Y(\ID)$ when
$|\cL|$ is the kind of linear system involved in Theorems (1.1) and
(1.2).
 \medbreak

{\bf Proposition (3.2)}\enspace {\it Set $\mu:=\mu(\ID)$.  Assume that
$\cL= \cM^{\ox m}\ox \cN$ where $\cN$ is spanned, $\cM$ is very ample,
and $m\ge \mu-1$.  If $Y(\ID)$ is nonempty, then it is locally closed in
$Y$, and is smooth and equidimensional.  Furthermore,}
        $$\cod(Y(\ID),Y)=\cod(\ID).$$ \medbreak

Indeed, let $y$ be an arbitrary point of $Y(\ID)$, and $C$ the
corresponding curve.  Consider the multigerm of $C$ along its singular
locus, and a corresponding miniversal deformation base space $B$.
(The formal analytic theory of $B$ was initiated by Schlessinger
\cite{Sc68}, and was algebraized by Artin \cite{At74} and Elkik
\cite{E73}.  An alternative complex analytic theory was developed by
Mather, Grauert, and others around the same time, and is explained in
\cite{T76, \S4}.)  The space $B$ is smooth, and its tangent space at the
origin is equal to
    $$H^0(\cO_S/\cK_{C,S}) \where \cK_{C,S}:=\cJ_{C,S}+\cI_{C,S},$$
 where $\cI_{C,S}$ is the ideal of $C$ on $S$.

There exists a map $\varphi$ from the germ of $Y$ at $y$ into $B$
carrying $y$ to the origin.  Since $m\ge\mu-1$ by hypothesis, $\varphi$ is
smooth, as we'll now show.  It suffices to show the surjectivity of the
map of tangent spaces, which is equal to the natural map,
    $$H^0(\cL)/H^0(\cO_S)\to H^0\bigl(\cL\big/\cK_{C,S}\cL\bigr).$$
 Since $\cN$ is spanned, it suffices to show the surjectivity of the map,
  $$H^0(\cM^{\ox m})\to
  H^0\bigl(\cM^{\ox m}\!\big/\cK_{C,S}\cM^{\ox m}\bigr).\eqno(3.2.1)$$

To show the surjectivity of (3.2.1), embed $S$ in a projective space $P$
so that $\cM=\cO_S(1)$, and let $\cK_{C,P}\subset\cO_P$  be the ideal
such that
  $$\cO_P/\cK_{C,P}=\cO_S/\cK_{C,S}.$$
  Now,
since $\cK_{C,S}\supset\cJ_{C,S}$, Proposition (3.1) implies that
        $\mu\ge\dim H^0(\cO_S/\cK_{C,S})$.
 Hence, by Gotzmann's regularity theorem \cite{Gtz78} (see also
\cite{Grn88, p.~80}), the ideal $\cK_{C,P}$ is $\mu$-regular.  So
$H^1(\cK_{C,P}(m))$ vanishes for $m\ge \mu-1$.  Hence the map
	$$H^0(\cO_P(m))\to
         H^0\bigl((\cO_P/\cK_{C,P})(m)\bigr)$$
 is surjective.  It factors through (3.2.1), so (3.2.1) is surjective.
Thus $\varphi$ is smooth.

The map $\varphi$ extends to a smooth map $U\to B$ where $U$ is a
complex analytic (or an \'etale) neighborhood of $y$ in $Y$.  The set
$U\cap Y(\ID)$ is carried onto the subset $B\es$ of $B$ whose points
represent the equisingular deformations of the multigerm of $C$.  (The
formal analytic theory of $B\es$ was initiated by Wahl \cite{W74}, and
can be algebraized using the work of Artin and Elkik cited above.  An
alternative complex analytic theory was developed by Teissier around the
same time, and is explained in \cite{T76, \S5}; this theory depends in
part on Wahl's, but not on Artin's and Elkik's.)  The subset $B\es$ is
closed in $B$, and is smooth.  Hence $Y(\ID)$ is locally closed in $Y$,
and is smooth.  Furthermore,
        $$\cod(Y(\ID),Y)=\cod(B\es,B).$$

Before proceeding with the proof of Proposition (3.2), note that,
together, the equation above and that in (3.2) yield the following
corollary, which gives another interpretation of $\cod(\ID)$.
 \medbreak

{\bf Corollary (3.3)}\enspace {\it Under the conditions of Proposition
\(3.2), the following formula gives the codimension of the subspace
$B\es$ of equisingular deformations in a miniversal deformation space
$B$ of the multigerm of any curve $C$ in $Y(\ID)$}:
        $$\cod(B\es,B)=\cod(\ID).$$ \medbreak

Returning to the proof of Proposition (3.2), recall that the analytic
multigerm of $C$ can be realized as the multigerm of a plane curve $C'$
of any suitably large degree $m'$.  (Many authors have written about
this assertion; for example, see \cite{Sh97} and its references.)  Let
$Y'$ be the projective space parameterizing all the plane curves of
degree $m'$, and let $y'\in Y'$ represent $C'$.  Shustin \cite{Sh91,
Thm.} proved that, locally at $y'$, the subset $Y'(\ID)$ is locally
closed in $Y'$, and is smooth; moreover, he gave a combinatorial
expression for $\cod(Y'(\ID),Y')$.  Applied to $Y',y'$ in place of
$Y,y$, the discussion before the corollary yields an alternative proof
that $B\es$ is closed in $B$, and is smooth.  (Shustin \cite{Sh91,
Rmk. 2} said that this fact about $B\es$ ``easily follows'' from his
theorem; doubtless, he envisioned a proof much like ours.)  It remains
to see that Shustin's combinatorial expression yields $\cod(\ID)$.

Shustin works, in fact, with plane curve germs, not multigerms; however,
the general case is an immediate consequence since the $B$, resp.\
$B\es$, of a multigerm is just the product of those of its component
germs.  Given a germ with central point $z$, Shustin \cite{Sh91, Def.~3,
p.~31} defines an integer $c(z)$ as follows.  If $z$ is simple, he sets
$c(z):=-2$.  Otherwise, he denotes the multiplicity of $z$ by $m$.  He
blows up the ambient plane at $z$, forms the exceptional divisor $E$
(which he denotes by $\Pi$), and lets $z_1,\dots,z_n$ be the points on
$E$ of the strict transform of the curve germ.  For each $i$, he lets
$d(z_i,E)$ be the number of points, including $z_i$, that are infinitely
near to $z_i$ and common to $E$ and the strict transform.
Finally, he sets
 $$c(z):=\tsum_{i=1}^n\bigl(c(z_i)+d(z_i,E)\bigr)+{m+1\choose 2}+n-3.$$

Shustin \cite{Sh91, Thm.} proved that $\cod(Y'(\ID),Y')=c(z)+1$.  So we
must prove that $c(z)+1=\cod(\ID)$, assuming $z$ is multiple.   Rewrite
the defining equation of  $c(z)$ as follows:
$$c(z)+1=\tsum_{i=1}^n\bigl(c(z_i)+1+d(z_i,E)\bigr)+{m+1\choose 2}-2.$$
 As in the proof of Proposition (3.1), say that $z_1,\dots,z_k$ are all
the simple points of the strict transform; furthermore, say that at
$z_1,\dots,z_l$ the branch is transverse to $E$, and at $z_{l+1},\dots,
z_k$ the branch is not transverse.  Finally, for $k+1\le i\le n$, let
$\ID_i$ denote the diagram associated to the germ at $z_i$ of the strict
transform.

If $1\le i\le l$, then $c(z_i)=-2$ and $d(z_i,E)=1$.  Hence
        $$c(z_i)+1+d(z_i,E)=0 \hbox{ if }1\le i\le l.$$
 If $z$ is an ordinary multiple point, then $l=n$; hence,
        $$c(z)+1=\textstyle{m+1\choose 2}-2=\cod(\ID).$$
 Proceeding by induction on the number of vertices in $\ID$, assume that
        $$c(z_i)+1=\cod(\ID_i)\hbox{ if }k+1\le i\le n.$$
 If $l+1\le i\le k$, then $c(z_i)=-2$.  Hence
 $$c(z)+1=\tsum_{i=l+1}^k\bigl(d(z_i,E)-1\bigr)
  +\sum_{i=k+1}^n\bigl(\cod(\ID_i)+d(z_i,E)\bigr)
        +{m+1\choose 2}-2.$$

Suppose $z_i$ is one of $z_{l+1},\dots, z_n$.  Then $z_i$ corresponds to
a vertex $S_i$ in $\ID$, and $S_i$ is an immediate successor of the
root, $R$ say.  Moreover, $d(z_i,E)$ is equal to the number of
successors $V_{ij}$ of $S_i$ that are proximate to $R$, and among the
$V_{ij}$, only $S_i$ is free.  If $l+1\le i\le k$, then $m_{V_{ij}}=1$
by Lemma (3.1).  If $k+1\le i\le n$, then $V_{ij}$ corresponds to a free
vertex in $\ID_i$, and $S_i$ corresponds to its root.  It follows that,
in the preceding display, the right hand side is equal to $\cod(\ID)$.
Thus $c(z)+1=\cod(\ID)$, as required.  We have now proved (3.2) and
(3.3). \medbreak

Set $d:=\deg(\ID)$.  This fifth character provides a sufficent condition
for $Y(\ID)$ to be nonempty; see Proposition (3.5) below.  To prove it,
our first step is to construct a weighted configuration $\cC$ of
infinitely near points of $S$ with $\ID$ as its associated diagram.

The construction is straightforward.  For every root of $\ID$, pick a
different point of $S$, and give it the corresponding weight.  Blowup
$S$ at each of these points.  On each exceptional divisor, pick a
different point for every immediate successor in $\ID$ of the
corresponding root, and give this infinitely near point the
corresponding weight.  Blow up at each of these points, and on each
exceptional divisor, pick a different point for every immediate
successor.  However, this time take care; if the successor is (remotely)
proximate to a root, then the corresponding new point must be taken on
the strict transform of the exceptional divisor belonging to this root.
Weight each new point appropriately, and continue the construction
similarly, taking suitable care with the satellites of $\ID$, until all
the vertices of $\ID$ have been handled.  This procedure is possible
because $\ID$ satisfies the Law of Succession.

Next, we pass from the configuration $\cC$ to an ideal $\cI$ on $S$,
which is {\it complete}, or integrally closed.  In the construction of
$\cC$, we perform a succession of blowups, each centered at a point
corresponding to a vertex of $\ID$.  In the end, we obtain a smooth
surface $S^*$.  Let $\beta\:S^*\to S$ be the natural map, and for each
vertex $V$ of $\ID$, let $E_V^*$ be the full inverse image on $S^*$ of
the point corresponding to $V$.  Finally, set
        $$\cI:=\beta_*\cO(-\tsum_{V\in\ID}m_VE_V^*).$$
 It is not hard to see that $\cI$ is complete.

The theory of complete ideals grew out of the study of complete linear
systems with base conditions on surfaces; see \cite{C97, \S8.3} for
example.  The theory contains a number of interesting and nontrivial
results, and two of them will be useful to us.  The first  asserts
the formula,
        $$\dim H^0(\cO_S/\cI)=d \where d:=\deg(\ID).\eqno(3.4)$$
 This formula provides an interpretation of the fifth character $d$.
The formula is a modern version of an old formula \cite{EC15, Vol. II,
\p.426}, and was proved independently, in various more general settings,
by Hoskin \cite{Ho56, 5.2, p.~85}, Deligne \cite{D73, 2.13, p.~22}, and
Casas \cite{C90}.  (Some authors attribute it only to Hoskin and
Deligne.)  It has been reproved and generalized by others.

The second result speaks to the existence, in the complete system
$|\cL|$, of a curve $C$ with $\cC$ as its configuration, so with $\ID$
as its diagram.  Denote by $Y(\cC)$ the subset of $Y$ parameterizing
these curves, and by $Y(\cI)$ the linear subspace of $Y$ associated to
$H^0(\cI\cL)$.  The result is this:\smallskip
 {\advance\leftskip by 2\parindent\noindent\it
 If $\cI\cL$ is generated by its global sections, then $Y(\cC)$ is a
{\rm nonempty} open subset of $Y(\cI)$. \smallskip}
 This result can be proved as follows.  On $S^*$, the image of
$\beta^*(\cI\cL)$ in $\beta^*\cL$ is an invertible sheaf, which is
generated by $H^0(\cI\cL)$.  So the corresponding linear system has no
base points.  Consider the smooth members that are transverse to the
components of the exceptional divisors $E_V$.  By Bertini's theorem,
these members form a nonempty open subset of $Y(\cI)$.  On the other
hand, it is not hard to see that this subset is equal to $Y(\cC)$.  (The
result is often stated in a local form, where $S$ is replaced by the
spectrum of a local ring; see \cite{C97, \S7.2} for example.)

Using the preceding result, we can prove the following proposition,
which gives a sufficent condition for $Y(\ID)$ to be nonempty.  The
proposition also gives another interpretation of $d$, which can be
paraphrased as follows: $d$ is the number of conditions imposed on the
members of $|\cL|$ by the requirement that they pass through the
infinitely near points of the configuration $\cC$ with multiplicities at
least as great as the corresponding weights. \medbreak

 {\bf Proposition (3.5)}\enspace {\it Assume that $\cL= \cM^{\ox m}\ox
\cN$ where $\cN$ is spanned, $\cM$ is very ample, and $m\ge d$ where
$d:=\deg(\ID)$.  Then $Y(\ID)$ is nonempty.  In fact, for any
configuration $\cC$ and ideal $\cI$ with diagram $\ID$ as above,
$Y(\cC)$ is a {\rm nonempty} open subset of $Y(\cI)$, and
$\cod(Y(\cI),Y)=d$.} \medbreak

Indeed, embed $S$ in a projective space $P$ so that $\cM=\cO_S(1)$, and
denote the preimage of $\cI$ in $\cO_P$ by $\cI_P$.  Then
$\cO_P/\cI_P=\cO_S/\cI$.  So $H^0(\cO_P/\cI_P)=d$ by Equation (3.4).
Hence $\cI_P$ is $m$-regular by Gotzmann's regularity theorem
\cite{Gtz78} (see also \cite{Grn88}, p.~80).  So $\cI_P(m)$ is generated
by $H^0(\cI_P(m))$, and $H^1(\cI_P(m))$ vanishes.
 
It follows that $\cI\cM^{\ox m}$ is generated by $H^0(\cI\cM^{\ox m})$,
and that the map,
        $$H^0(\cM^{\ox m})\to H^0(\cM^{\ox m}\!\big/\cI\cM^{\ox m}),$$
 is surjective.

Since $\cN$ is spanned, it follows that $\cI\cL$ is generated by
$H^0(\cI\cL)$ and that the sequence,
        $$0\to H^0(\cI\cL)\to H^0(\cL)\to H^0(\cL\big/\cI\cL)\to0,$$
 is exact, also on the right.  This exactness implies that $Y(\cI)$ has
codimension $d$ in $Y$, and the result displayed just above implies that
$Y(\cC)$ is a nonempty open subset of $Y(\cI)$.  The proof is now
complete.

We turn now to the sixth and last character $\dim(\ID)$.  Although $\ID$
is fixed, we make a number of choices when we construct $\cC$.  In fact,
we choose a point in $S$ for each root of $\ID$, and a point in a
certain exceptional divisor for each remaining free vertex.  The choices
are not completely arbitrary.  However, it is intuitively clear that the
number of ``degrees of freedom'' is just twice the number of roots plus
the number of remaining free vertices, or $\dim(\ID)$.

To formalize the preceding discussion, note that the correspondence,
$\cC\mapsto \cI$, is injective.  Indeed, $Y(\cC)$ is a {\rm nonempty}
open subset of $Y(\cI)$ for a suitable $\cL$ thanks to Proposition
(3.5), or thanks simply to the displayed result preceding it.
Conversely, every complete ideal arises from some configuration,
although its minimal Enriques diagram is not necessarily $\ID$; see
\cite{C97, \S8.3} for example.

Form the Hilbert scheme $\Hilb^d_S$.  In it, there is a point for every
ideal $\cI$, thanks to Equation (3.4).  The various points form a
subset, and it is locally closed by virtue of the work of Nobile and
Villamayor \cite{NV97, Thm.~2.6} or that of Lossen [pvt.\ comm.].
Denote it by $H(\ID)$, and view $H(\ID)$ as a reduced subscheme.  In
terms of $H(\ID)$, we have the following precise interpretation of the
sixth and last character of $\ID$.  \medbreak

 {\bf Proposition (3.6)}\enspace {\it In $\Hilb^d_S$, the subscheme
$H(\ID)$ is smooth and equidimensional.  Furthermore,
        $\dim H(\ID)=\dim(\ID)$.
 } \medbreak

The proposition is intuitively clear from the discussion preceding it.
This discussion is developed into a formal proof in \cite{K--P}, but the
proof is surprisingly long and involved.  (However, it also shows that
$H(\ID)$ is irreducible if $S$ is; moreover, it works in arbitrary
characteristic.)  Alternatively, the proposition can be derived from
Proposition (3.2), Proposition (3.5), and the following lemma; we'll
prove the lemma, and then derive the proposition. \medbreak

{\bf Lemma (3.7)}\enspace {\it Let $Z(\ID)$ in $Y\x H(\ID)$ denote the
incidence scheme, and $Z_0(\ID)$ the preimage of $Y(\ID)$ in $Z(\ID)$.
Then $Z_0(\ID)\to Y(\ID)$ is bijective.

Assume that $\cL= \cM^{\ox m}\ox \cN$ where $\cN$ is spanned and $\cM$
is very ample.  Set $\mu:=\mu(\ID)$ and $d:=\deg(\ID)$.  If $m\ge d$,
then $Z(\ID)/H(\ID)$ is a bundle of projective spaces.  If
$m\ge\max(d,\,\mu-1)$, then $Z_0(\ID)$ is a dense open subset of
$Z(\ID)$, and the projection $Z(\ID)\to Y$ induces an isomorphism
$Z_0(\ID)\risom Y(\ID)$.}\medbreak

Indeed, given $y\in Y(\ID)$, let $C$ be the corresponding curve, $\cC$
its configuration, $\cI$ the complete ideal of $\cC$, and $h\in H(\ID)$
the point representing $\cI$.  Then $(y,h)\in Z(\ID)$.  So $Z_0(\ID)\to
Y(\ID)$ is surjective.  Now, given $(y,h')\in Z(\ID)$, let $\cI'$ be the
complete ideal represented by $h'$, and $\cC'$ the corresponding
configuration.  Every point $P'$ in $\cC'$ lies on a strict transform
$C'$ of $C$, and at $P'$ the multiplicity of $C'$ is at least the weight
of $P'$ in $\cC'$.  It follows that $\cI'$ contains $\cI$.  However,
both these ideals have colength $d$ in $\cO_S$.  Hence the two ideals
are equal.  Hence $h'=h$.  Thus, $Z_0(\ID)\to Y(\ID)$ is bijective.

By general principles, for any $\cL$, there is a coherent sheaf $\cQ$ on
$H(\ID)$ such that $Z(\ID)$ is equal to $\IP(\cQ)$.  In our case, if
$m\ge d$, then the fibers of $Z(\ID)/H(\ID)$ are all of the same
dimension by Proposition (3.5).  Hence $\cQ$ is locally free because, by
hypothesis, $H(\ID)$ is reduced.  Thus $Z(\ID)/H(\ID)$ is a bundle of
projective spaces.

Assume $m\ge\max(d,\,\mu-1)$.  Then $Z_0(\ID)$ is locally closed as
$Y(\ID)$ is so by (3.2).  On the other hand, every fiber of $Z(\ID)\to
H(\ID)$ is a projective space, and meets $Z_0(\ID)$ in a nonempty open
subset by (3.5).  Hence $Z_0(\ID)$ is dense in $Z(\ID)$.  Therefore,
since $Z_0(\ID)$ is locally closed, it is open.  Now, the map
$Z_0(\ID)\to Y(\ID)$ is birational since it is bijective (and the
characteristic is zero).  Hence, by Zariski's Main Theorem, $Z_0(\ID)\to
Y(\ID)$ is an isomorphism since it is bijective and since $Y(\ID)$ is
smooth by (3.2).  Thus the lemma is proved.

To derive Proposition (3.6), apply Lemma (3.7) with $m\ge\max(\mu-1,d)$
and $\cL= \cM^{\ox m}$ where $\cM$ is very ample.  Then $Z_0(\ID)$ is
smooth and equidimensional since $Y(\ID)$ is so by (3.2).  Moreover,
$Z_0(\ID)\to H(\ID)$ is smooth and surjective, since (3.5) and (3.7)
imply that $Z(\ID)\to H(\ID)$ is so, since every fiber meets $Z_0(\ID)$,
and since $Z_0(\ID)$ is open.  Hence, $H(\ID)$ is smooth and
equidimensional.  Finally, $\dim H(\ID)$ is equal to the difference
between $\dim Z_0(\ID)$ and the dimension of the fibers of $Z_0(\ID)\to
H(\ID)$.  Hence the equation $\dim H(\ID)=\dim(\ID)$ follows from the
equations in (3.2) and (3.5).  Thus (3.6) is proved, and all six
characters of $\ID$ have been interpreted.

\sct4 Proofs of the theorems

In this section we prove Theorems (1.1) and (1.2).  Our proofs use a
recursive formula in $r$ for the cycle class representing the $r$-nodal
curves in an algebraic system on an algebraic family of surfaces.  Even
though our ultimate interest lies in a linear system on a fixed surface,
the added generality is a necessity, not a luxury.  Indeed, we proceed
inductively by passing to a new system on a family of blowups of the
initial family of surfaces.  The new system is not linear, and the new
family is not constant, even when the initial system and family are so.

Let $\pi\:F\to Y$ be a smooth and projective family of (possibly
reducible) surfaces, where $Y$ is equidimensional and Cohen--Macaulay,
and let $D$ be a relative effective divisor on $F/Y$.  Denote by $p_j:
F\x_Y F\to F$ the $j$th projection, by $\Delta \subset F\x_Y F$ the
diagonal subscheme, and by $\cI_\Delta$ its ideal.  Then $D$ is defined
by a global section $s$ of the invertible sheaf $\cO_F(D)$, and $s$
induces a section $s_i$ of the sheaf of relative twisted principal
parts,
        $$\cP_{F/Y}^{i-1}(D)
        :=p_{2*}\bigr(p_1^*\cO_F(D)\big/\cI_\Delta^i\bigr)\for i\ge1.$$
Denote the scheme of zeros of $s_i$ by $X_i$.  So
$X_1=D$.  Furthermore, as a set, $X_i$ consists of the points $x\in F$
at which $D_{\pi(x)}$ has multiplicity at least $i$.  As $i$ varies, the
$X_i$ form a descending chain of closed subschemes.

The sheaf $\cP_{F/Y}^{i-1}(D)$ fits into the exact sequence,
        $$0\to\cSym^{i-1}\Omega^1_{F/Y}(D)\to\cP_{F/Y}^{i-1}(D)
        \to\cP_{F/Y}^{i-2}(D)\to0,$$
where the first term is the symmetric power of the sheaf of relative
differentials.
 Hence $\cP_{F/Y}^{i-1}(D)$ is locally free of rank $i+1\choose 2$ by
induction on $i$.  Therefore, every component of $X_i$ has codimension
at most $i+1\choose 2$.  Furthermore, if every component has exactly
this codimension, then the fundamental class $[X_i]$ is equal to the top
Chern class of $\cP_{F/Y}^{i-1}(D)$, and so $[X_i]$ can be expressed as
a polynomial in the following three Chern classes:
  $$v:=c_1(\cO_F(D))\and w_j:=c_j(\Omega^1_{F/Y})\for j=1,2.\eqno(4.1)$$

 For example, in this way we obtain the following three expressions,
which we need for the proofs of Theorems (1.1) and (1.2):
 $$\setbox0=\hbox{$:=$}\def\&{\cr&\kern\wd0} \eightpoint\eqalign{
 [X_2]&=v^3+w_1v^2+vw_2;\cr
 [X_3]&=v^6+4w_1v^5+(5w_1^2+5w_2)v^4+(2w_1^3+11w_1w_2)v^3
        +(6w_2w_1^2+4w_2^2)v^2+4vw_1w_2^2;\cr
 [X_4]&=v^{10}+10w_1v^9+(15w_2+40w_1^2)v^8+(82w_1^3+111w_1w_2)v^7\&
         +(91w_1^4+315w_2w_1^2
        +63w_2^2)v^6+(52w_1^5+429w_2w_1^3+324w_1w_2^2)v^5\&
         +(12w_1^6+282w_2w_1^4+593w_2^2w_1^2+85w_2^3)v^4
         +(72w_2w_1^5+464w_2^2w_1^3+259w_1w_2^3)v^3\&
        +(132w_2^2w_1^4
        +246w_2^3w_1^2+36w_2^4)v^2
         +(72w_2^3w_1^3+36w_1w_2^4)v.\cr}$$

Let $b\:F'\to F\x_Y F$ be the blowup along $\Delta$, and set $p_j'
:= p_j\circ b$.  Then $p_2'\: F'\to F$ is again a smooth and
projective family of surfaces; in fact, over a point $x$ of $F$, the
fiber $F'_x:= p_2^{\prime-1}x$ is just the blowup of the fiber $F_{\pi x}:=
\pi^{-1}\pi x$ at $x$.  Set $F_i:=p_2^{\prime-1}(X_i)$, and let $\pi_i:F_i\to
X_i$ be the restriction of $p_2'$ to $F_i$.  In sum, we have the
following diagram:
        $$\xymatrix{
  F\ar[d]_\pi & F\x_Y F\ar[l]_{p_1}\ar[d]^{p_2} & {F'}
        \ar[l]_{\quad b}\ar[d]_{p_2'}
         & F_i\ar@{_{(}->}[l] \ar[d]_{\pi_i} \\
  Y & F\ar[l]_\pi & F\ar@{=}[l] & X_i\ar@{_{(}->}[l]
        }$$

The pullback $p_1^{-1}D$ is defined by $p_1^*s$.  The restriction of
$p_1^*s$ to $p_2^{-1}X_i$ is not simply a section of the restriction of
$p_1^*\cO_F(D)$, but is also a section of the restriction of its
subsheaf $\cI_\Delta^i(p_1^{-1}D)$, in view of the definition of $X_i$.
Hence the restriction of $p_1^{\prime *}s$ to $F_i$ is a section of the
restriction of $\cO_{F'}(p_1^{\prime-1}D-i b^{-1}\Delta)$.  Therefore, the
difference of the restrictions,
        $$D_i:=(p_1^{\prime-1}D)|F_i-(i b^{-1}\Delta)|F_i,$$
 is a relative effective divisor on $F_i/X_i$.

Let $\ID$ be an abstract minimial Enriques diagram.  In Section 3, we
associated some loci to it, and made a study of them.  We'll now
continue that study in the present context, and apply it to prove
Theorems (1.1) and (1.2).  First, set $d:=\deg(\ID)$, and in the
relative Hilbert scheme $\Hilb^d_{F/Y}$, form the subset $H_{F/Y}(\ID)$
parameterizing the complete ideals sitting on the fibers of $F/Y$ and
having $\ID$ as diagram.  Just as before, $H_{F/Y}(\ID)$ is locally
closed, and we'll view it as a reduced subscheme.  Next, in $Y$ form the
subset $Y(\ID)$ of points over which the fiber of $D$ is a reduced curve
with $\ID$ as diagram.  Next, form the scheme
        $$Z(\ID):=H_{F/Y}(\ID)\textstyle\bigcap\Hilb^d_{D/Y},$$
 and form the set-theoretic preimage $Z_0(\ID)$ of $Y(\ID)$ in $Z(\ID)$.
Finally, form the set-theoretic image of $Z(\ID)$ in $Y$, and the
image's closure $U(\ID)$.

Assume for a moment that $F=S\x Y$ where $S$ is a smooth, irreducible,
projective surface.  Assume that $Y$ is the parameter projective space
of a complete linear system $|\cL|$ on $S$, and that $D\subset F$ is the
total space.  Then $\Hilb^d_{F/Y}$ is equal to $Y\x\Hilb^d_S$; so
$H_{F/Y}(\ID)$ is equal to $Y\x H(\ID)$ where $H(\ID)$ is the locus of
Section 3.  Furthermore, $Y(\ID)$, $Z(\ID)$, and $Z_0(\ID)$ are equal to
their counterparts in Section 3, but $U(\ID)$ is new.  If we assume in
addition that $\cL$ is suitably ample, then we may say more about these
loci, as the following proposition asserts.
 \medskip
 {\bf Proposition (4.2)}\enspace {\it Preserve the preceding conditions.
Let $1\le r\le8$.  Assume that $\cL= \cM^{\ox m}\ox \cN$ where $\cN$ is
spanned, $\cM$ is very ample, and $m\ge 3r$.  Set $s:=\cod(\ID)$.
 \part1 Assume $s\le r$.  Then $Y(\ID)$ is locally closed, smooth, and
equidimensional, and $\cod(Y(\ID),Y)=s$.  Also, $Z(\ID)$ is reduced,
$Z_0(\ID)$ is dense and open in $Z(\ID)$, and the projection induces an
isomorphism of reduced schemes, $Z_0(\ID)\risom Y(\ID)$.  Finally,
$Y(\ID)$ is open and dense in $U(\ID)$.
 \part2 Independent of $\ID$, there is a closed subset $Y'$ of $Y$ with
$\cod(Y',Y)\ge r+1$ such that $Y'$ contains $Y(\ID)$ if $s\ge r+1$.
Moreover, $Y'$ contains every point of $Y$ over which the fiber of $D$
is a curve with a multiple component.}
 \medskip
 Indeed, let's first check the following two inequalities when $s\le10$:
        $$\mu(\ID)-1\le s \and \deg(\ID)\le 3s.\eqno(4.3)$$
 Both sides of both inequalities are ``additive,'' so we may assume that
$\ID$ has only one root.  Then both inequalities can be checked by
inspecting Table 2-1.  (In fact, the second inequality holds without any
bound on $s$; see \cite{K--P}.)

Combining Display (4.3) with Propositions (3.2) and (3.5) and with Lemma
(3.7), we obtain the first two assertions in Part (1).  So $Z_0(\ID)$
and $Z(\ID)$ have the same closure in $Y\x\Hilb^d_S$.  This closure
projects onto $U(\ID)$ since $\Hilb^d_S$ is projective.  Now, $Z_0(\ID)$
maps onto $Y(\ID)$.  Hence $Y(\ID)$ is dense in $U(\ID)$.  Therefore,
since $Y(\ID)$ is locally closed, it is open in $U(\ID)$.

To prove Part (2), note that there are only finitely many minimal
Enriques diagrams $\ID'$ with $\deg(\ID')\le 3r$.  For each such $\ID'$,
form $Z(\ID')$ in $Y\x H(\ID')$, and $U(\ID')$ in $Y$.  Then
$Z(\ID')/H(\ID')$ is a fiber bundle by Lemma (3.7).  The dimensions of
its fibers and of its base are given by Propositions (3.5) and (3.6).
Hence,
        $$\cod(U(\ID'),Y)\ge\cod(\ID').$$
 If $\cod(\ID')\le r$, then equality holds, and $Y(\ID')$ is open and
dense in $U(\ID')$ by Part (1) applied with $\ID'$ for $\ID$.

Let $Y''$ be the union of the $U(\ID')$ with $\cod(\ID')\ge r+1$.
Let $Y'''$ be the union of the differences $U(\ID')-Y(\ID')$ with
$\cod(\ID')= r$.  Set $Y':=Y''\cup Y'''$.  Then $\cod(Y',Y)\ge r+1$.
To prove that $Y'$ contains $Y(\ID)$ if $s\ge r+1$, we need to develop a
little more of the general theory of diagrams.

Call a minimal Enriques diagram $\ID'$ a {\it subdiagram\/} of $\ID$,
and say that $\ID$ {\it contains\/} $\ID'$, if $\ID'$ consists of some
of the same vertices, equipped with possibly smaller weights and with
the induced relations of succession and proximity.  For example,
$A_{2i-1}$ is a subdiagram of $A_{2i}$ and of $A_{2i+1}$, but $A_{2i}$
is not a subdiagram of $A_{2i+1}$.  Also, $D_{2i}$ is a subdiagram of
$D_{2i+1}$ and of $D_{2i+2}$, but $D_{2i+1}$ is not a subdiagram of
$D_{2i+2}$.  The following lemma asserts the existence of a convenient
subdiagram (it is proved in \cite{K--P} for any $r\ge1$).
 \medskip
{\bf Lemma (4.4)}\enspace {\it Let $1\le r\le8$.  If $\cod(\ID)\ge r+1$,
then $\ID$ contains a minimal Enriques subdiagram $\ID'$ such that
$\cod(\ID')\ge r$ and $\deg(\ID')\le 3r$.}\medskip

Assuming the lemma, let's complete the proof of the proposition.  Assume
$s\ge r+1$.  Let $y$ be an arbitrary point of $Y(\ID)$, and $C$ the
corresponding curve.  Use the notation and constructions of Sections 2
and 3.  Form the corresponding weighted configuration $\cC$ of
infinitely near points of $S$.  The lemma provides a certain subdiagram
$\ID'$ of $\ID$.  Let $\cC'$ be the corresponding subconfiguration of
$\cC$.  Let $\cI$ and $\cI'$ be the associated complete ideals on $S$.
Letting $m'_V$ denote the weight of a vertex $V\in\ID'$ we have
$$\cI:=\beta_*\cO\bigr(-\tsum_{V\in\ID}m_VE_V^*\bigl)
        \subset\beta_*\cO\bigr(-\tsum_{V\in\ID'}m'_VE_V^*\bigl)=\cI';$$
 the latter equality does not hold by definition, but can be derived
using the Leray Spectral Sequence.

The inclusion $\cI\subset\cI'$ induces an inclusion,
        $H^0(\cI\cL)\subset H^0(\cI'\cL)$.
 In turn, the latter yields  $Y(\cI)\subset Y(\cI')$.
However, $y\in Y(\cI)$ and $Y(\cI')\subset U(\ID')$.  So, if
$\cod(\ID')\ge r+1$, then $y\in Y''$.  Suppose $\cod(\ID')= r$.  Then
$\ID'\neq \ID$ since $s\ge r+1$.  Hence $y\notin Y(\ID')$.  So $y\in
Y'''$.  Since $y$ is arbitrary, therefore $Y(\ID)\subset Y'$.

Finally, let $y$ be a point of $Y$ over which the fiber of $D$ is a
curve with a multiple component.  On that component, pick $r$ distinct
points, and let $\cI'$ be the product of the squares of their maximal
ideals.  Then $\cI'$ is a complete ideal with diagram $\ID':=rA_1$.
Now, $\cod(\ID')= r$ and $\deg(\ID')= 3r$.  Moreover, $y\in Y(\cI')$ and
$y\notin Y(\ID')$.  So $y\in Y'''$.  Thus Proposition (4.2) is proved,
given Lemma (4.4).

To prove Lemma (4.4), note that $\ID$ always contains $A_1$, and recall
from Table 2-1 that $\deg(A_1)=3$.  Hence, if $r=1$, we may take $A_1$
as $\ID'$.  So assume $r\ge2$.  Proceed by induction on the number of
roots.

First, assume that $\ID$ has only one root, say $R$ of weight
$m$.  If $m\ge 4$, then we may take $\ID'$ to be $X_{1,0}$ if $r\ge4$
since $\cod X_{1,0}=8$ and $\deg X_{1,0}=10$; moreover, then we may take
$\ID'$ to be $D_4$ if $r=3,2$ since $\cod D_4=4$ and $\deg D_4=6$.  If
$m=2$, then $\ID$ is $A_k$ with $k\ge r+1$, and we may take $\ID'$ to be
$A_r$ if $r$ is odd and to be $A_{r+1}$ if $r$ is even.

Suppose $m=3$.  Then $\ID$ is of type $D$, $E$, or $J$.  If $\ID$ is
either $J_{l,j}$ or $E_{6l+j}$ where $l\ge2$, then $\ID$ begins with a
succession of two vertices of weight 3.  Hence, then we may take $\ID'$
to be $J_{2,0}$ if $r\ge4$ since $\cod J_{2,0}=10$ and $\deg
J_{2,0}=12$.  Moreover, then we may take $\ID'$ to be $A_3$ if $r=3,2$
since $\cod A_3=3$ and $\deg A_3=6$.

Suppose that $\ID$ is either $E_7$ or $E_8$.  Then $r\le7$.  Also, $E_8$
contains $E_7$, and $E_7$ contains $A_3$.  Hence, we may take $\ID'$ to
be $E_7$ if $r\ge4$ since $\cod E_7=7$ and $\deg E_7=10$.  Moreover, we
may take $\ID'$ to be $A_3$ if $r=3,2$.

Suppose $\ID$ is $E_6$.  Then $r\le5$.  Also, $E_6$ contains $A_2$.
Hence, we may take $\ID'$ to be $E_6$ if $r\ge3$ since $\cod E_6=6$ and
$\deg E_6=9$.  Moreover, we may take $\ID'$ to be $A_2$ if $r=2$ since
$\cod A_2=2$ and $\deg A_2=5$.

Finally, assume that $\ID$ has more than one root.  Let $\ID_1$ be a
connected component, and $\ID_2$ its complement; so
        $$\ID= \ID_1+ \ID_2.$$
 Set $s_i:=\cod(\ID_i)$.  Now, $\ID_1$ has only one root.  If $s_1\ge
r+1$, then by the one-root case, already $\ID_1$ contains a $\ID'$ such
that $\cod(\ID')\ge r$ and $\deg(\ID')\le 3r$.  If $s_1=r$, then
$\deg(\ID_1)\le 3s_1$ by the second inequality of Display (4.3) applied
to $\ID_1$; so we may take $\ID_1$ as $\ID'$.

Suppose $s_1\le r-1$.  Set $r_2:=r-s_1$.  Then $s_2\ge r_2+1$ because
        $$s_1+s_2=\cod(\ID_1)+\cod(\ID_2)=\cod(\ID)\ge r+1,$$
 by the definition of $s_i$, by the additivity of `cod', and by
hypothesis.  So, by induction, we may assume that $\ID_2$ contains a
$\ID'_2$ such that $\cod(\ID'_2)\ge r_2$ and $\deg(\ID'_2)\le 3r_2$.  Set
        $$\ID':=\ID_1+\ID_2'.$$
 Then by the additivity of `cod' and of `deg' we have
 $$\eqalign{\cod(\ID')&{}=\cod(\ID_1)+\cod(\ID_2')\ge s_1+r_2=r, \and\cr
         \deg(\ID')&{}=\deg(\ID_1)+\deg(\ID_2')\le 3s_1+3r_2=3r,\cr}$$
 the last inequalities coming from Display (4.3).  The lemma is now
proved.

Let $\pi\:F\to Y$ and $D$ be arbitrary again.  Let $1\le r\le8$.  In
Proposition (4.2), the hypotheses concerning $\cL$ no longer apply, but
in both Parts (1) and (2), the assertions still make sense.  Make the
genericity hypothesis on $F/Y$ and $D$ that these assertions hold.

The genericity hypothesis implies immediately that the class,
        $$u(D,s):=[U(sA_1)] \for 1\le s\le r,$$
 is of codimension $s$ on $Y$, and that it enumerates the $s$-nodal
curves in the family $D/Y$, since they form a dense subset of $U(sA_1)$.
For convenience, extend the definition of $u(D,s)$ by setting
        $$u(D,0):=1, \and u(D,s):=0 \for s\le -1.$$
 Our next step is to express $u(D,r)$ in terms of the three Chern
classes $v$, $w_1$ and $w_2$ of Display (4.1).  We'll find that $u(D,r)$
can be expressed as a polynomial of degree $r$ in the following classes:
   $$y(a,b,c):=\pi_*v^aw_1^bw_2^c \where a,b,c\ge0 \and a+b+2c\le r+2.
        \eqno(4.5)$$
 We'll just sketch the procedure here; more details are found in
\cite{K--P}.

Consider the closed subset $Y'$ of $Y$ described in Part (2) of
Proposition (4.2).  Since $\cod(Y',Y)\ge r+1$, we may replace $Y$ by
$Y-Y'$ without affecting the validity of our expression for $u(D,r)$.
Thus we may (and will) assume that every fiber of $D/Y$ is reduced and
its diagram has $\cod$ at most $r$.

The genericity hypothesis also implies that $X_2$ is reduced,
Cohen--Macaulay, and equidimensional of codimension 3 in $F$.  Indeed,
there is a natural embedding of $F$ in $\Hilb^3_{F/Y}$, given by sending
a point of $F$ to the square of its maximal ideal, and the image of $F$
is simply $H_{F/Y}(A_1)$.  Hence, $X_2$ is isomorphic to $Z(A_1)$, so
reduced of codimension 3.  Similarly, $X_3$ is reduced of codimension 6
if $4\le r\le8$, and $X_4$ is reduced of codimension 10 if $r=8$.  Of
course, $X_3$ is empty if $r\le3$, and $X_4$ is empty if $r\le7$ since
the fibers of $D/Y$ are now well behaved.

The genericity hypothesis further implies that the analogous genericity
hypothesis holds for $F_i/X_i$ and $D_i$ for $i=2,3$ provided $r$ is
replaced by $r-1$, respectively $r-4$, and this number is at least 1.
However, there is one exception: over every point of $X_4$, the fiber of
$D_2$ contains twice the exceptional divisor, and so is nonreduced.
Moreover, we have the following recursive formula:
  $$ru(D,r) = \pi_*\bigl(u(D_2,r-1)-u(D_3,r-4)+3281u(D_4,r-8)\bigr).
        \eqno(4.6)$$
 All these statements are discussed in detail in \cite{K--P}; here,
we'll accept them as true.

Intuitively, the preceding statements make sense.  For example, on the
left side of  Formula (4.6), the $r$ is present because, over a point of
$Y(rA_1)$, the fiber of $X_2$ consists of $r$ distinct points.  On the
right, the term $u(D_2,r-1)$ enumerates the $(r-1)$-nodal curves
$C$ in the family $D_2/X_2$.  If $C$ is transverse to the
exceptional curve $ E$ in its ambient fiber of $F_2/X_2$, then
$C$ is the proper transform of an $r$-nodal curve in the family
$D/Y$, the $r$th node being at the blowup center.  If a $C$ were
tangent to $ E$, then it would be the transform of a curve with a
cusp and $(r-1)$ nodes; however, the locus of these curves is
$Y(A_2+(r-1)A_1)$, so empty since $Y'$ is now empty.  Similarly, no
$C$ has one node on $ E$ and $(r-1)$ nodes off it.

However, a $C$ can have three nodes on $ E$ if $4\le r\le8$.
Then $C=  E+ C'$ where $C'$ is an $(r-4)$-nodal
curve transverse to $ E$.  Such a $C'$ is the proper
transform of a curve with an ordinary triple point at the blowup center.
So $C'$ is an $(r-4)$-nodal curve in the family $D_3/X_3$.  These
curves are enumerated by the class $u(D_3,r-4)$.  Thus $u(D_3,r-4)$
appears with a minus sign on the right in Formula (4.6).

The multiplier 3281 is more problematical.  Certainly, when $r=8$, some
multiple of $u(D_4,r-8)$ will appear in Formula (4.6), but the
multiplier might well vary with $F/Y$ and $D$.  In \cite{K--P}, residual
intersection theory is used to prove that the multiplier does not vary;
then its value is determined from a particular case.

To express $u(D,r)$ as a polynomial in the classes $y(a,b,c)$ of Display
(4.5), we proceed formally by recursion on $r$ using Formula (4.6) and
the formulas for the  $[X_i]$ displayed early in the section.
First,  for $r=1$, we obtain
        $$\eqalign{
        u(D,1)&=\pi_*u(D_2,0)=\pi_*[X_2]=\pi_*(v^3+w_1v^2+vw_2)\cr
        &=y(3,0,0)+y(2,1,0)+y(1,0,1).\cr}$$

To proceed, set $e:=[b^{-1}\Delta]$, which is the class of the
exceptional divisor.  To lighten the notation, let $e$, $v$, and $w_j$
also denote their own pullbacks.  Then
        $$c_1(\cO_{F_i}(D_i))=v-ie,\ c_1(\Omega^1_{F'/F})=w_1+e,
        \and c_2(\Omega^1_{F'/F})=w_2-e^2.$$
 So, using Formula (4.6) and recursion in $r$, express $u(D,r)$ as a
polynomial of degree $r-1$ in terms of the form,
        $$\pi_*(\pi_i)_*((v-ie)^a(w_1+e)^b(w_2-e^2)^c).$$
 Expand each term, and reduce the powers of $e$ using the basic relation,
        $$e^3+w_1e^2+w_2e=0.$$
 Push out to $F'$, and use the following three identities:
        $$ b_*[F_i]=[F\x X_i]\and  b_*e=0\and  b_*e^2=-[\Delta].$$
 Thus $u(D,r)$ can be expressed as a polynomial of degree $r$ in the
$y(a,b,c)$.

The actual polynomial expression can be found much more efficiently and
written much more compactly as follows.  As in Section 1, define
polynomials $P_r(a_1,\dots,a_r)$ by the formal identity in $t$,
 $$\textstyle\sum_{r\ge0}P_rt^r/r!
        = \exp\bigl(\sum_{q\ge1}a_q t^q/q!\bigr).$$
 This time, however, the variables $a_q$ will not be replaced by
numbers, but by classes on $Y$.  These classes are certain linear
combinations of the $y(a,b,c)$ of Display (4.5), and are computed as
follows.

Given a formal monomial $(v-ie)^a(w_1+e)^b(w_2-e^2)^c$, expand it.  Reduce
the result modulo the relation $e^3=-w_1e^2-w_2e$, collect the terms
involving $e^2$, and denote the coefficient of $e^2$ by $q_i(a,b,c)$.  Set
        $$Qx_i(v^aw_1^bw_2^c)=-q_i(a,b,c).$$
 Extend the definition of $Qx_i(\bullet)$ by
linearity to all polynomials in $v$, $w_1$, $w_2$.

Define polynomials $b_q$ in $v$, $w_1$, $w_2$ recursively by setting
$b_1:=[X_2]$ and
 $$\eqalign{b_{q+1}
        :={}&P_q(Qx_2(b_1),\dots,Qx_2(b_q)) [X_2] \textstyle
  -3!\,{q\choose 3}P_{q-3}(Qx_3(b_1),\dots,Qx_3(b_{q-3})) [X_3]\cr
        &\qquad+3281\cdot7!\,\textstyle{q\choose 7}
         P_{q-7}(Qx_4(b_1),\dots,Qx_4(b_{q-7})) [X_4],\cr}$$
 where the $[X_i]$ are the polynomials in $v$, $w_1$, $w_2$ introduced
early in the section.  Finally, take $v$, $w_1$ and $w_2$ to be the
Chern classes defined in Display (4.1), and set
        $$a_q:=\pi_*(b_q) \for q=1,\dots,r.$$
 Then it can be checked by brute force that
        $$u(D,r)=P_r(a_1\dots,a_r)/r!\,.$$
 It would be good to have a conceptual proof that our two mechanical
procedures always produce the same polynomials.

We are now ready to prove Theorems (1.1) and (1.2).  So take $F=S\x Y$
where $S$ is the given surface and where $Y$ is the $n$-dimensional
parameter projective space of the given complete linear system $|\cL|$
on $S$.  Take $D\subset F$ to be the total space of $|\cL|$.  Then
$\cO_F(D)$ is equal to $\cL\ox\cO_Y(1)$, and so
        $$v=l+h \where l:=c_1(\cL) \and h:=c_1(\cO_Y(1)).$$
 Also, $\Omega^1_{F/Y}$ is simply the pullback of $\Omega^1_S$.  Hence,
only four of the classes $y(a,b,c)$ are nonzero, namely,
        $$y(r+2,0,0),\ y(r+1,1,0),\ y(r,2,0),\ y(r,0,1),$$
 and their values are, respectively,
        $$\textstyle{r+2\choose2}dh^r,\ (r+1)kh^r,\ sh^r,\ xh^r,$$
 where $d$, $k$, $s$, and $x$ are the four numbers defined at the
beginning of Section 1.

Using the procedure of the preceding paragraph, we get
        $$u(D,r)=(P_r(d,k,s,x)/r!)h^r$$
 where $P_r$ is viewed as a polynomial in $d,k,s,x$ as in Theorem (1.1).
Finally, $N_r$ is just the degree of $u(D,r)$ because of Proposition
(4.2) and of the following little lemma; apply the lemma with $Y(rA_1)$
as $U$.  Thus Theorem (1.1) is proved.

 \medskip {\bf Lemma (4.7)}\enspace {\it In the above setup, let $U$ be
a reduced subscheme of $Y$ of pure codimension $r$, and $U'$ its
boundary.  Assume that $U$ parameterizes curves without multiple
components, and set $m:=n-r$.  Then there exists a nonempty open subset
$S_m$ of $S^{\x m}$ with this property: let $M$ be the linear $r$-space
representing those curves that pass through the points of any given
$m$-tuple in $S_m$; then $M\cap U$ is finite and reduced, and $M\cap U'$
is empty.}\medskip

 Indeed, in $D$, form the smooth locus $D_0$ of the projection $D\to Y$.
Form the fibered product $D_0^{\x m}$.  Then $D_0^{\x m}\to Y$ is smooth.
Hence $D_0^{\x m}\x_YU$ is reduced.  Moreover, it is dense in $D^{\x
m}\x_YU$.  Also, $D^{\x m}\x_YU$ is of pure dimension $2m$.  Consider
the natural map,
        $$D^{\x m}\x_YU\to S^{\x m}.$$
 By Sard's theorem, its fibers are finite and reduced over a nonempty
open subset $S_m$ of $S^{\x m}$.  These fibers are simply the $M\cap U$.

On the other hand, $D^{\x m}\x_YU'$ is of dimension at most $2m-1$.  So
the map,
         $$D^{\x m}\x_YU'\to S^{\x m},$$
 cannot be surjective.  Hence, its fibers are empty over a nonempty
subset of $S_m$.  Replace $S_m$ by this subset.  Then these fibers are
the $M\cap U'$.  Thus the lemma is proved.

To prove Theorem (1.2), we work with the  following seven classes on $Y$,
        $$[U(D_4+(r-4)A_1)]\for r=4,5,6,7,\
         [U(D_6)],\ [U(D_6+A_1)],\ [U(E_7)].$$
 Again, Proposition (4.2) and Lemma (4.7) imply that the degrees of
these classes are just the numbers we seek.

{}From the discussion in the second paragraph after (4.6), we get the
equality,
        $$[U(D_4+(r-4)A_1)]=u(D_3,r-4).$$
 That discussion is set-theoretic, but the corresponding cycles are
reduced because the requisite genericity hypothesis holds; it holds
initially because of Proposition (4.2), and subsequently because
genericity propagates.  Now, the class $u(D_3,r-4)$ was computed
implicitly above as part of the recursion that led to formulas in
Theorem (1.1).  By making it explicit for $r=4,5,6,7$, we obtain the
first four formulas in Theorem (1.2).

Looking back, we see that we've established an equation of the form,
 $$\textstyle\sum_{r\ge0}u(D_3,r)t^r/r!
        = \exp\bigl(\sum_{q\ge1}a_q t^q/q!\bigr),$$
 at least modulo $t^4$, where the $a_q$ are linear combinations of
certain classes on $X_3$, namely, of the analogues for $F_3/X_3$ and
$D_3$ of the $y(a,b,c)$ of Display (4.5).

We obtain the last three formulas in Theorem (1.2) from the following
three formulas, which we'll prove in a moment:
        $$\eqalignno{
     [U(D_6)]&=\pi_*[U_2(D_4+A_1)]-2[U(D_4+2A_1], &(4.8.1)\cr
        [U(D_6+A_1)]&=\pi_*[U_2(D_4+2A_1)]-3[U(D_4+3A_1], &(4.8.2)\cr
        [U(E_7)]&= \pi_*[U_2(D_6)]-[U(D_6+A_1)], &(4.8.3)\cr}$$
 where $U_2(\ID)$ is the analogue of $U(\ID)$ for $F_2/X_2$ and $D_2$.
We just discussed how to find the class $[U(D_4+jA_1)]$ for $j=0,1,2,3$.
By applying that result with $F/X$ and $D$ replaced by $F_2/X_2$ and
$D_2$, we obtain an expression for $[U_2(D_4+jA_1)]$.  Then the first
two formulas above yields expressions for $[U(D_6)]$ and $[U(D_6+A_1)]$.
Applying the first of these with $F/X$ and $D$ replaced by $F_2/X_2$ and
$D_2$, we obtain an expression for $[U_2(D_6)]$.  Then the third formula
above yields an expression for $[U(E_7)]$.

It remains to prove Formulas (4.8.1), (4.8.2) and (4.8.3).  Let $C$ be a
curve (or fiber) in the family $D_2/X_2$ with diagram $D_4+(j+1)A_1$
where $0\le j\le1$.  If $C$ does not contain the exceptional curve $E$
in its ambient fiber of $F_2/X_2$, then $C$ is the proper transform of a
curve that is in the family $D/Y$ and that has, at the blowup center, a
singularity with diagram $A_k$ where $1\le k\le3$.  Therefore, the full
diagram of this curve is one of the following:
        $$D_4+(j+2)A_1,\ D_4+(j+1)A_1+A_2,\hbox{ or }D_4+jA_1+A_3.$$
 On the other hand, if $C$ contains $E$, then since $C$ is reduced,
$C=C'+E$ where $C'$ is the proper transform of a curve that has, at the
blowup center, a singularity with diagram $D_6$.

Thus $X_2(D_4+(j+1)A_1)$ is carried by $\pi$ into the union of the sets
$Y(\ID)$ where $\ID$ is one of the following four diagrams:
 $$D_4+(j+2)A_1,\ D_4+(j+1)A_1+A_2,\ D_4+jA_1+A_3\hbox{ or }D_6+jA_1.$$
 In fact, $X_2(D_4+(j+1)A_1)$ is carried onto the union because the
preceding analysis is reversible.  Of these $\ID$, only $D_4+(j+2)A_1$
and $D_6+jA_1$ have $\cod(\ID)=6+j$; the other two $\ID$ have
$\cod(\ID)\ge7+j$.  Moreover, above a point in $Y(D_4+(j+2)A_1)$ there
are $j+2$ points in $X_2(D_4+(j+1)A_1)$, whereas above a point in
$Y(D_6+jA_1)$, there is only one point.  It now follows from Proposition
(4.2) that Formulas (4.8.1) and (4.8.2) hold.

Formula (4.8.3) may be proved similarly.  Let $C$ be a curve in the
family $D_2/X_2$ with diagram $D_6$.  If $C$ does not contain $E$, then
$C$ is the proper transform of a curve that is in $D/Y$ and that has, at
the blowup center, a singularity with diagram $A_k$ where $1\le k\le3$.
If $C$ contains $E$, then $C=C'+E$ where $C'$ is the proper transform of
a curve that has, at the blowup center, a singularity with diagram
$E_7$.  In fact, this analysis is reversible.  Thus $X_2(D_6)$ is
carried by $\pi$ onto the union of the sets $Y(\ID)$ where $\ID$ is one
of the following three diagrams:
        $$D_6+A_1,\ D_6+A_2,\hbox{ or } E_7. $$
 Now, $\cod(D_6+A_1)=7$ and $\cod(E_7)=7$, but $\cod(D_6+A_2)=8$.
Moreover, above a point in either $Y(D_6+A_1)$ and $Y(E_7)$, there is
only one point in $X_2(D_6)$.  It now follows from Proposition (4.2)
that the final formula (4.8.3) holds.

Thus both Theorems (1.1) and (1.2) are proved!

\sct5 Vainsencher's treatment

Fix $r\le7$.  Vainsencher \cite{V95} enumerated the $r$-nodal curves in
a ``suitably general'' linear system of dimension $r$ on a smooth,
irreducible, projective surface $S$.  By ``suitably general,'' he meant
that the system is a subsystem of a complete system $|\cL|$, and as such
corresponds to a point in a suitable nonemtpy open subset of the
Grassmannian of subsystems; in addition, $\cL=\cM^{\ox m}$ where $\cM$
is ample and $m\ge m_0$ for a suitable $m_0$.  He established the
existence of the open subset and of $m_0$ if $r\le6$.  If $r=7$, he said
in his Section 7 that then his argument does not apply.  Moreover, he
said that his final formula does not yield integers in certain examples;
however, he later found a computational error [pvt. comm., Dec. 1997].

In this section, we relate Vainsencher's enumeration to our own.
Thereby, we justify his approach on the basis of ours for every $r$,
including $r=7$.  Furthermore, we prove that the term ``suitably
general'' may be used in the following more precise sense, which was
specified in Theorem (1.1) of the introduction: the system is {\it suitably
general\/} if it is a subsystem defined by the condition to pass through
$n-r$ general points where $n$ is the dimension of $|\cL|$, and if
$\cL=\cM^{\ox m}\ox \cN$ where $\cN$ is spanned, $\cM$ is very ample,
and $m\ge 3r$.  Of course, it follows formally from our result that
there exists a nonemtpy open subset of the Grassmannian of subsystems of
$|\cL|$, whose points represent systems that work.

Vainsencher used an ad hoc recursive definition of the type of a
singularity of a curve $C$ on $S$.  Namely, he called a sequence of
points $(x_1,x_2,\dots,x_r)$ a ``singularity of type''
$(m_1,m_2,\dots,m_r)$ of $C$ if $x_1$ is a point of $C$ of multiplicity
at least $m_1$ and if, on the blowup of $S$ at $x_1$, the sequence
$(x_2,\dots,x_r)$ is a singularity of type $(m_2,\dots,m_r)$ on the
difference $C_1-m_1E_1$, where $C_1$ is the total transform of $C$, and
$E_1$ is the exceptional divisor.  To avoid confusion with the
traditional notion of singularity, call $(x_1,x_2,\dots,x_r)$ a {\it
singularity sequence}.

Vainsencher said that $(x_1,x_2,\dots,x_r)$ is of {\it strict type\/} if
$x_1$ is a point of $C$ of multiplicity exactly $m_1$, if $x_2$ lies off
of $E_1$, and if $(x_2,\dots,x_r)$ is of strict type on $C_1-m_1E_1$;
otherwise, he said that the sequence is of {\it weak type}.  He wrote
$m^{[s]}$ in place of a string of $s$ repetitions of $m$; for example,
$(3,2^{[3]})$ stands for $(3,2,2,2)$.

Fix a sequence $(x_1,\dots,x_r)$ of type $(2^{[r]})$.  Distinguish the
$x_i$ that actually lie in $C$, and at their union, consider the
singularity of $C$.  Let $\ID$ be its minimal Enriques diagram, and set
$k:=\cod(\ID)$.  We are going to prove that $k\ge r$, and to determine
when $k=r$.  To do so, assume that $k\le r$; we are going to analyze the
cases, and to prove that $k=r$ or to find a contradiction.  To begin,
assume in addition that only one of the $x_i$ lies in $C$; necessarily,
it is $x_1$.

Suppose $x_1\in C$ is a singularity of type $A_k$.  Then $x_1$ is
resolved by a sequence of either $k/2$ or $(k+1)/2$ blowups at centers
that are points of multiplicity 2.  Hence $(k+1)/2\ge r$.  However,
$k\le r$.  Hence $k=r=1$.  Thus $x_1$ is of type $A_1$, and $(x_1)$ is
of strict type (2).

Suppose $x_1\in C$ is of type $D_k$ or $E_k$; these are the only other
types of singularities with $k\le7$.  Then $x_1$ is of multiplicity 3.
So $C_1-2E_1$ is equal to $C'+E_1$ where $C'$ is the strict transform of
$C$.  We now prove that $x_1$ must be of type $D_4$, $D_6$, or $E_7$.

Suppose $x_1$ is of type $D_k$.  Then $C'+E_1$ has, along $E_1$, a
singularity of type $3A_1$ if $k=4$, of type $A_1+A_3$ if $k=5$, of type
$A_1+D_4$ if $k=6$, or of type $A_1+D_5$ if $k=7$.  Hence, if $k=4$,
then $r\le4$, and so $k=r$ since $k\le r$.  Moreover, it follows that,
if $k=6$, then $r\le6$, and so $k=r$.  Now, an $A_3$ singularity is
resolved by two blowups at centers that are points of multiplicity 2.
Hence, if $k=5$, then $r\le4$, a contradiction.  Moreover, it follows
that, if $k=7$, then $r\le6$, a contradiction.

Suppose $x_1$ is of type $E_k$ with $k=6$ or $k=7$.  Then $C'+E_1$ has a
singularity at $x_2$ of type $A_5$ if $k=6$, or of type $D_6$ if $k=7$;
moreover, $x_3,\dots,x_r$ lie infinitely near $x_2$.  Hence, it follows
from the analysis above of types $A_5$ and $D_6$ that $k=r=7$.

Suppose now that several of the $x_i$ lie in $C$.  For each $x_i$ in
$C$, let $\ID_i$ be the diagram of the singularity of $C$ at $x_i$, set
$k_i=\cod(\ID_i)$, and denote by $r_i$ the number of $x_j$ that lie
infinitely near $x_i$, including $x_i$.  Apply the above analysis to the
sequence of these $x_j$.  It yields $k_i\ge r_i$, and enumerates the
possible types when $k_i=r_i$.  Now, $\sum k_i=k$ and $\sum r_i=r$.
Hence, $k\ge r$, and if $k=r$, then $k_i=r_i$.  Putting it all together,
we obtain the second column of Table 5-1, which lists the possible types
of singularities when $k=r$.

$$\matrix\multispan3\hfil\bf Table \number\sctno-1\hfil\cr
\multispan3\hfil\bf
 Sequences of type $(2^{[r]})$ with $\cod(\ID)=r$ \hfil\strut\cr
\noalign{\smallskip\hrule\smallskip}
 r &\hbox{type of singularity}&\hbox{\# of permutations}\strut\cr
\noalign{\smallskip\hrule\smallskip}
 1,\ 2,\ 3 & rA_1 & r! \cr
 4 & 4A_1,\ D_4 & 4!,\ 6\cr
 5 & 5A_1,\ D_4+A_1 & 5!,\ 30\cr
 6 & 6A_1,\ D_4+2A_1,\ D_6 &6!,\ 180,\ 30\cr
 7 & 7A_1,\ D_4+3A_1,\ D_6+A_1,\ E_7 &7!,\ 1260,\ 210,\
 30\cr\endmatrix
 $$
 \bigskip
Given a singularity of a type listed in Table 5-1, there is a sequence
$(x_1,\dots,x_r)$ of type $(2^{[r]})$ giving rise to it.  Indeed, a
second look at the analysis above shows the existence of a sequence.
The set of points is uniquely determined, but their order is
subject to permutation.  The number of permutations is not hard to find;
it is listed, case by case, in the third column of the table.  (In
\cite{L98}, 
 Liu provided a combinatorial device, ``admissible graphs,'' to help
find them.) Two examples are worked out next to illustrate how to find
these numbers directly.

For example, a $D_4$ singularity $x_1\in C$ arises from the sequence
$(x_1,\dots,x_4)$ where $x_2$, $x_3$, and $x_4$ are the three points
infinitely near $x_1$.  The three points are subject to 6 permutations.
Similarly, a $D_4+A_1$ singularity arises from the sequence
$(x_1,\dots,x_5)$ where $x_1,\dots,x_4$ are as before and $x_5$ is a
node.  The four points $x_2,\dots,x_5$ are nodes on the blowup of $S$ at
$x_1$; so the four are subject to 24 permutations.  In addition, the
five points $x_1,\dots,x_5$ may be cyclically permuted to
$x_5,x_1,\dots,x_4$, and then the three points $x_2,x_3,x_4$ are subject
to 6 permutations.  In total, the sequence $(x_1,\dots,x_5)$ is subject
to 30 permutations.

Let $C$ vary now in a linear system that is suitably general in the
sense specified above in the second paragraph of the section.  Then, by
Section 4, only finitely many of these $C$ have a diagram $\ID$ with
$\cod(\ID)=r$, and no $C$ has $\cod(\ID)>r$.  Hence there is a finite
number of singularity sequences of type $(2^{[r]})$ on all the $C$;
denote this number by $N(2^{[r]})$.  Therefore, with the notation of
Section 1, we find the following formulas:
 $$\eqalign{N(2^{[r]})&=r!\,N_r \for r=1,2,3,\cr
 N(2^{[4]})&=4!\,N_4+6\,N(3),\cr
 N(2^{[5]})&=5!\,N_5+30\,N(3,2),\cr
 N(2^{[6]})&=6!\,N_6+180\,N(3,2,2)+30\,N(3(2)),\cr
 N(2^{[7]})&=7!\,N_7+1260\,N(3,2,2,2)+210\,N(3(2),2)+30\,N(3(2)').\cr}
        \eqno(5.1)$$

Vainsencher found similar formulas; see Proposition 4.1 and Section 7 in
\cite{V95}.  However, instead of the coefficients 180 and 1260, he had
90 and 210 (or literally, 1260/6), since he enumerated, not the curves
with singularities of types $D_4+2A_1$ and $D_4+3A_1$, but the sequences
of strict types $(3,2,2)$ and $(3,2,2,2)$, which are subject to $2!$ and
$3!$ permuations respectively.

Consider the set of all singularity sequences of type $(2^{[r]})$.  It
underlies a natural scheme, which is formed, using the setup of Section
4, as follows.  If $r=1$, take the scheme $X_2$.  If $r=2$, take the
corresponding scheme of the family $D_2/X_2$.  If $r\ge3$, repeat this
procedure $r-2$ more times.  Now, using the methods of Section 4, it is
not hard to find an expression in terms of the four numbers $d,k,s,x$ of
Section 1 for the degree of this scheme, and to check that this
expression is equal to the polynomial obtained from the above formulas
and from Theorems (1.1) and (1.2).  It follows in particular that the scheme is
reduced, and that its degree (when considered as a cycle on $Y$) is equal
to $N(2^{[r]})$; the reducedness
also follows directly from Section 4.

Vainsencher obtained his formulas for the $N_r$ from his version of the
formulas in Display (5.1) and his version of Theorem (1.2).  Thus we
have indeed justified and extended his approach.

The formulas in Display (5.1) are implicit in our proof of Theorem (1.1)
in Section 4, and it is interesting to make them explicit as we are
about to do.  In particular, we recover the coefficients using less
combinatorics.  Furthermore, as in our proof of Theorem (1.1), we must
work in a more general setup: we must let $C$ vary in an algebraic
system on a family of surfaces satisfying the genericity hypothesis that
the assertions in Parts (1) and (2) of Proposition (4.2) both hold.
Thus, in this more general setup, we now establish the corresponding
versions of the formulas in Display (5.1).

First of all, if $r\le4$, then the formula for $N(2^{[r]})$ results
immediately from the recursive formula (4.6) by taking degrees.

To obtain the formula for $N(2^{[5]})$ in Display (5.1), apply the
formula for $N(2^{[4]})$ to $D_2/X_2$, getting
        $$N(D_2;2^{[4]})=4!\,N_4(D_2)+6\,N(D_2;3).$$
  Identifying the terms, we find
        $$N(2^{[5]})=4!\deg\pi_*u(D_2,4)+6\,N(3,2).$$
 Hence (4.6) yields
$$N(2^{[5]})=\bigl(5!\,N_5+24\,N(3,2)\bigr)+6\,N(3,2)=5!\,N_5+30\,N(3,2),$$
 which is the formula in Display (5.1).

Similarly, to obtain the formula for $N(2^{[6]})$, apply the formula for
$N(2^{[5]})$ to $D_2/X_2$, getting
        $$N(D_2;2^{[5]})=5!\,N_5(D_2)+30\,N(D_2;3,2).$$
 Now, $N(D_2;3,2)=2\,N(3,2,2)+\,N(3(2))$ follows from Formula (4.8.1) by
taking degrees.  Identifying the several terms, we find
        $$N(2^{[6]})=5!\,\deg\pi_*u(D_2,5)+60\,N(3,2,2)+30\,N(3(2)).$$
 Hence (4.6) yields
$$N(2^{[6]})=\bigl(6!\,N_6+120\,N(3,2,2)\bigr)+60\,N(3,2,2)+30\,N(3(2)).$$
 Combining the second and third terms yields the formula in Display
(5.1).

Finally, to obtain the formula for $N(2^{[7]})$, apply the formula for
$N(2^{[6]})$ to $D_2/X_2$, and identify the terms, getting similarly
 $$\eqalign{N(2^{[7]})&=\bigl(7!\,N_7+720\,N(3,2,2,2\bigr))
        +180\bigl(3\,N(3,2,2,2)+N(3(2),2)\bigr)\cr
        &\qquad+30\bigl(N(3(2),2)+N(3(2)')\bigr).\cr}$$
 Combining terms yields the formula in Display (5.1), and completes our
work.

\references
 \serial{ajm}{Amer. J. Math.}
 \serial{cmp}{Comm. Math. Phys.}
 \serial{duke}{Duke Math. J.}
 \serial{ens}{Ann. Sci. \'Ecole Norm. Sup.}
  \serial{ja}{J. Alg.}
  \serial{imrn}{Int. Math. Res. Notices}
 \serial{invm}{Invent. Math.}
 \serial{jag}{J. Alg. Geom.}
  \serial{ja}{J. Alg.}
  \serial{jams}{J. Amer. Math. Soc.}
 \serial{jram}{J. Reine Angew. Math.}
 \serial{ma}{Math. Ann.}
 \serial{mz}{Math. Z.}
 \serial{plms}{Proc. London Math. Soc.}
 \serial{sms}{Selecta Math. Sovietica}
 \serial{sdg}{Surveys in Diff. Geom.}
 \serial{tams}{Trans. AMS}

A76
 V. I. Arnold,
 Local normal forms of functions,
 \invm 35 1976 87--109.

At74
 M. Artin,
 Versal deformations and algebraic stacks,
 \invm 27 1974 165--189.

BL97
 J. Bryan and C. Leung,
 The enumerative geometry of K3 surfaces and modular forms,
 alg-geom/9711031 $=$ \jams 13(2) 2000 371--410.

BL98
 J. Bryan and C. Leung,
 Generating functions for the number of curves on Abelian surfaces,
 math.AG/9802125 $=$ \duke 99{\rm (2)} 1999 311--28.

BL98s
 J. Bryan and C. Leung,
 Counting curves on irrational surfaces,
 \sdg 5 1999 313-39.

CH96
 L. Caporaso and J. Harris,
 Counting plane curves of any genus,
 \invm 131 1998{\rm (2)} 345--92.

C90
 E. Casas,
 Infinitely near imposed singularities and singularities of polar
curves, \ma 287 1990 429--54.

C97
 E. Casas,
 Singularities of plane curves,
 London Math Society Lecture Note Series {\bf276}, Cambridge Univ. Press, 2000.

D73
 P. Deligne,
 Exp. X, Intersections sur les surfaces r\'eguli\`eres,
 {\it in} ``SGA 7 II, Groupes de Monodromie en g\'eom\'etrie
alg\'ebrique,''
 Lecture Notes in Math. {\bf 340}, Springer-Verlag, 1973, 1--38.

E73
 R. Elkik,
 Solutions d'\'equations \`a co\'efficients dans un anneau hens\'elien,
\ens 6 1973 553--604.

EC15
 F. Enriques and O. Chisini,
 Lezioni sulla teoria geometrica delle equazioni e delle funzioni
algebriche, Zanichelli, Bologna, 1915.

Gtt98
 L. G\"ottsche,
 A conjectural generating function for numbers of curves on surfaces,
 \cmp 196 1998 523--533.

Gtz78
 G. Gotzmann,
 Eine Bedingung f\"{u}r die Flachheit und das Hilbertpolynom eines
gradu\-ierten Ringes,
 \mz 158 1978{\rm (1)} 61--70.

Grn88
 M. Green,
 Restrictions of linear series to hyperplanes, and some results of
Macaulay and Gotzmann,
 {\it in} ``Algebraic Curves and Projective Geometry (Trento 1988),''
(E. Ballico and C. Ciliberto, eds.), Lecture Notes in Math. {\bf 1389},
Springer-Verlag, 1989, 76--86.

GLS98
 G.-M. Greuel, C. Lossen, and E. Shustin,
 Plane curves of minimal degree with prescribed singularities,
 \invm 133 1998 539--580.

Hi53
 F. Hirzebruch,
 \"Ubertragung einiger S\"atze aus der Theorie der algebraischen
Fl\"achen auf komplexe Mannigfaltigkeiten von zwei komplexen
Dimensionen,
 \jram 191 1953 110--124.

Ho56
 M. A. Hoskin,
 Zero-dimensional valuation ideals associated with plane curve branches,
 \plms 6 1956 70--99.

KPB
 S. Kleiman and R. Piene,
 Node polynomials: results and examples,
 to appear.

K--P
 S. Kleiman and R. Piene,
 Node polynomials for curves on surfaces,
 to appear.

L98
 A. K. Liu,
 Family blowup formula, admissible graphs and the counting of nodal
 curves, math.AG/9804095.

NV97
 A. Nobile and O. E. Villamayor,
 Equisingular stratifications associated to families of planar ideal,
 \ja 193 1997 239--59.

Ran89
 Z. Ran,
 Enumerative geometry of singular plane curves,
 \invm 97 1989 447--469.

Sc68
 M. Schlessinger,
 Functors of Artin rings,
 \tams 130 1968 208--222.

Sh91
 E. I. Shustin,
 On manifolds of singular algebraic curves,
 \sms 10 1991 27--37.

Sh97
 E. I. Shustin,
 Smoothness of equisingular families of plane algebraic curves
 \imrn 2 1997 67--82.

S76
 D. Siersma,
 Periodicity in Arnold's list of singularities,
 {\it in} ``Real and complex singularities, Nordic Summer School, Oslo
1976,'' P.  Holm (ed.), Sijthoff \& Noordhoff 1977, 87--109.

T76
 B. Teissier,
 The hunting of invariants in the geometry of discriminants, {\it in}
 ``Real and complex singularities, Nordic Summer School, Oslo 1976,''
 P.  Holm (ed.), Sijthoff \& Noordhoff 1977, 565--677.

V95
 I. Vainsencher,
 Enumeration of $n$-fold tangent hyperplanes to a surface,
 \jag 4 1995 503--526.

V97
 R. Vakil,
 Counting curves of any genus on rational ruled surfaces,
 alg-geom/9709003.

W74
 J. Wahl,
 Equisingular deformations of plane algebroid curves,
\tams  193 1974 143--170.

Z65
 O. Zariski,
 Studies in equisingularity I,
 \ajm 87 1965 507--536 \ = Coll. Papers {\bf IV}, 31--60.

\endreferences 
\enddocument